\documentclass[reqno,twoside,12pt]{amsart}
\usepackage{amssymb,amsfonts,amsthm,amsmath}
\usepackage{enumitem}
\usepackage[pagebackref=false]{hyperref}
\usepackage[hmargin=1in,vmargin=1in]{geometry}
\usepackage{cite}

\def\eqdef{\stackrel{\rm def}{=}}

\def\d{{\rm d}}
\def\ddt{\frac{\d}{\d t}}

\newcommand{\tnum}{\rm(\roman*)}
\newcommand{\rnum}{\rm(\alph*)}

\def\beq{\begin{equation}}
\def\eeq{\end{equation}}
\def\beqs{\begin{equation*}}
\def\eeqs{\end{equation*}}

\newtheorem{theorem}{Theorem}[section]
\newtheorem{lemma}[theorem]{Lemma}
\newtheorem{proposition}[theorem]{Proposition}
\newtheorem{corollary}[theorem]{Corollary}
\newtheorem{definition}[theorem]{Definition}
\newtheorem{assumption}[theorem]{Assumption}

\theoremstyle{definition}
\newtheorem{remark}[theorem]{Remark}
\newtheorem{example}[theorem]{Example}
\newtheorem*{xnotation}{Notation}

\usepackage{color}

\def\varep{\varepsilon}

\newcommand{\R}{\ensuremath{\mathbb R}}

\newcommand{\N}{\ensuremath{\mathbb N}}

\newcommand{\bigo}{\mathcal O}

\numberwithin{equation}{section}

\title{On the finite time blow-ups  for solutions of nonlinear differential equations}
\author{Luan Hoang}

\address{Department of Mathematics and Statistics,
Texas Tech University\\
1108 Memorial Circle, Lubbock, TX 79409--1042, U. S. A.}
\email{luan.hoang@ttu.edu}


\makeatletter
\@namedef{subjclassname@2020}{\textup{2020} Mathematics Subject Classification}
\makeatother

\date{\today}

\begin{document}

\begin{abstract} 
We study systems of nonlinear ordinary differential equations  where the dominant term, with respect to large spatial variables, causes blow-ups and is positively homogeneous of a degree $1+\alpha$ for some $\alpha>0$. We prove that the asymptotic behavior of a solution $y(t)$ near  a finite blow-up time $T_*$ is $(T_*-t)^{-1/\alpha}\xi_*$ for some nonzero vector $\xi_*$. Specific error estimates for  $|(T_*-t)^{1/\alpha}y(t)-\xi_*|$  are provided. In some typical cases, they can be a positive power of $(T_*-t)$ or $1/|\ln(T_*-t)|$. This depends on whether the decaying rate of the lower order term, relative to the size of the dominant term, is of a power or logarithmic form. Similar results are obtained for a class of nonlinear differential inequalities with finite time blow-up solutions. Our results cover larger classes of nonlinear equations, differential inequalities and error estimates than those in the previous work. 
\end{abstract}

\maketitle

\tableofcontents

\pagestyle{myheadings}\markboth{L. Hoang}
{On the finite time blow-ups for  solutions of nonlinear differential equations}

\section{Introduction}\label{intro}
In the article \cite{FS84a}, Foias and Saut study the Navier--Stokes equations (NSE) with potential body forces in the  functional form 
\beq\label{NSE}
u'+Au+B(u,u)=0
\eeq
in an infinite dimensional space. 
Above, $A$ is  the (linear) Stokes operator which has positive eigenvalues and $B(\cdot,\cdot)$ is a bilinear form. It is proved in \cite{FS84a} that any nontrivial solution  $u(t)$ of \eqref{NSE} has the following asymptotic behavior
\beq\label{FSlim}
e^{\Lambda t}u(t)\to \xi_* \text{ in any $C^m$-norms as $t\to\infty$,}
\eeq
where $\Lambda$  an eigenvalue of $A$ and  $\xi_*$ is an eigenfunction of $A$ associated with $\Lambda$.

This result and its proof in \cite{FS84a} are later extended to abstract differential inequalities in  Hilbert spaces in \cite{Ghidaglia1986a,Ghidaglia1986b}. 
We recall here the main result in \cite{Ghidaglia1986a}, but streamline and restrict it to $\R^n$ to fit the scope of the current paper.

Suppose $A$ is a self-adjoint $n\times n$ (real) matrix and $u:[0,\infty)\to \R^n$ satisfies 
\beq\label{Dineq}
|u'+Au|\le |u| \zeta(t) \text{ for all $t\in (0,\infty)$. }
\eeq

If $\zeta(t)$ is integrable and square integrable on $(0,\infty)$, then the limit in \eqref{FSlim} holds true in $\R^n$ for some eigenvector $\xi_*$ of $A$. 

The proof in \cite{FS84a}, and also in \cite{Ghidaglia1986a}, is based on the analysis of the asymptotic behaviors of the Dirichlet quotient and the normalized solution, see the corresponding definitions in \eqref{quot} below. Another method of proving \eqref{FSlim} for ordinary differential equations (ODEs) is presented in \cite{CaHK1}.
The result \eqref{FSlim} is the starting point for a theory of  asymptotic expansions for the solutions of the NSE in \cite{FS87}. See \cite{FS91,FGuillope86,FS84b,FS86,FHS1,FHOZ1,FHOZ2,FHN1,FHN2,HM1} for its further developments for the NSE including the associated nonlinear spectral manifolds and normal forms.   
These asymptotic expansions are also established for other general nonlinear ODEs  and partial differential equations (PDEs)  in \cite{Minea,Shi2000,HTi1,CaHK1}.  
See also \cite{CaH1,CaH2,H6,HM2} for results for the NSE with non-potential body forces, and \cite{CaH3,H5} for ODE systems with forcing functions.

In general, an asymptotic expansion can be obtained independently from the limit \eqref{FSlim}.
However, in the case the nonlinearity has positively homogeneous terms of  noninteger degrees such as in \cite{CaHK1},  \eqref{FSlim} is the essential first step and $\xi_*$ is critical in obtaining other terms in the expansion.

All equations in the papers cited above share the same feature, namely, they have  a \emph{linear} lowest order term of dissipative type. This linear term  dominates the other nonlinear terms for small solutions.
For ``genuine" nonlinear equations, the paper \cite{H7} studies a nonlinear ODE systems in $\R^n$ with  a superlinear lowest order term. More specifically, it considers the equation
\beq\label{oldeq}
y'=-H(y) Ay +G(t,y),
\eeq
where the matrix $A$ has positive eigenvalues,
the function $H$ is positive and positively homogeneous of a degree $\alpha>0$, see Definition \ref{phom} below, and the function $G$ is a higher order term. 
It is proved in \cite{H7} that any nontrivial, decaying solution of \eqref{oldeq} satisfies
\beqs
\lim_{t\to\infty} t^{1/\alpha} y(t) =\xi_* \text{ which is an eigenvector of $A$.}
\eeqs

In the recent work \cite{H8}, the function $H$ in \eqref{oldeq} is assumed to be positively homogeneous of a negative degree $(-\alpha)$. The behavior of a solution $y(t)$ near a  finite extinction time $T_*$ is proved to be 
\beqs
\lim_{t\to T_*^-} (T_*-t)^{-1/\alpha} y(t) =\xi_* \text{ which is an eigenvector of $A$.}
\eeqs

In both papers \cite{H7,H8}, the solution $y(t)$ goes to zero as time tends to infinity or a finite value.
In contrast, the current paper is devoted to studying the solutions that blow up in finite time.
We consider  the following ODE system in $\R^n$
\beq\label{yGeq}
y'=H(y) Ay +G(t,y),
\eeq
where $H$ is as in \cite{H7}, while $G(t,y)$ represents  a lower order term. 

Regarding the differential inequalities, we also study
\beq\label{newdi}
|y'-H(y) Ay|\le  |y|^{1+\alpha}E_0(t),\text{ for some nonnegative function $E_0(t)$.}
\eeq
This clearly is an an extension of the inequality \eqref{Dineq}.

We study the solutions $y(t)$ of \eqref{yGeq} or \eqref{newdi} that blow up at a finite time $T_*$.
We will prove that, under appropriate conditions,  there exists a nonzero vector $\xi_*$ such that
\beq\label{mainlim}
\lim_{t\to T_*^-} (T_*-t)^{1/\alpha}y(t)=\xi_*
\eeq
and will  provide explicit error estimates for $|(T_*-t)^{1/\alpha}y(t)-\xi_*|$. 
Some typical cases can be described roughly as follows. For equation \eqref{yGeq},
if
\beq\label{Gdec}
\frac{ |G(t,x)| }{ |x|^{1+\alpha} } = \bigo(|x|^{-\delta})
\text{ or } \bigo((\ln |x|)^{-\delta}) \text{ as $|x|\to\infty$,}
\eeq 
then
\beq\label{rough}
|(T_*-t)^{1/\alpha}y(t)-\xi_*|=\bigo((T_*-t)^{\delta'})\text{ or }\bigo(|\ln(T_*-t)|^{-\delta'})
\text{ as $t\to T_*^-$, respectively.}
\eeq
Similarly for inequality \eqref{newdi}, if 
\beq\label{etadec}
E_0(t) = \bigo((T_*-t)^\delta) \text{ or }\bigo(|\ln(T_*-t)|^{-p}) \text{ for $p>p_0$, as $t\to T_*^-$,  } 
\eeq
where $p_0$ is a certain number in $(1,\infty)$,
then one has \eqref{rough} again.

The presence of the logarithmic functions in \eqref{Gdec}--\eqref{etadec} is new compared to previous work \cite{H7,H8}.

For a convenient connection between \eqref{yGeq} and \eqref{newdi}, we convert inequality \eqref{newdi} to the equation
\beq\label{mainode}
y'=H(y) Ay + f(t)
\eeq
 with   $|f(t)|\le |y(t)|^{1+\alpha}E_0(t)$.
Therefore, \eqref{mainode} is also an extension of \eqref{Dineq}. Moreover, we will use \eqref{mainode} as a tool to analyze \eqref{yGeq}.

This paper is organized as follows.
In Section \ref{blowtime}, we impose the general conditions on the matrix $A$ and function $H$ in Assumption \ref{assumpAH}. The function $G(t,y)$ in \eqref{yGeq}  satisfies condition \eqref{Grate} which is more general than the usually assumed upper bound $|y|^{1+\alpha-\delta}$. We prove in Theorem \ref{blowthm} that any solution of \eqref{yGeq}  with sufficiently large initial data blows up in finite time. 
In Section \ref{simeg} we study equation \eqref{basicf} which will  serve  as a reference system for our perturbation technique in Section \ref{symcase} later.
We establish the precise asymptotic behavior \eqref{mainlim} for its blow-up solutions in Theorem \ref{simthm}.
 In Section \ref{ratesec}, we study  system \eqref{mainode} and  derive the lower and upper bounds for its solutions which are comparable with the exact blow-up rate $(T_*-t)^{-1/\alpha}$. These bounds are crucial to many estimates in the remainder of the paper.
Sections \ref{symcase} and \ref{nonsymcase} contain the main results for solutions of equation \eqref{mainode}, see Theorems \ref{symthm} and \ref{nonsymthm}. 
We establish the limit \eqref{mainlim} with the error estimates being specified in details in  \eqref{remy}--\eqref{yxio} and \eqref{symo1}--\eqref{symo3}.
Because of the nonlinear function $H$ and the general forcing function $f(t)$, many technical conditions are imposed, see, e.g., \eqref{condE}, Assumptions \ref{E2assum}  and \ref{moreE} and \eqref{tEfin}. 
Nonetheless, these technicalities will provide flexibility in obtaining different results in  Section \ref{egsec}.
The final Section \ref{egsec} applies the general results in the previous two sections to more specific cases of equations \eqref{mainode} and \eqref{yGeq}. 
Equation  \eqref{mainode} is dealt with in Part I in which two typical forms of $f(t)$ result in  Theorems \ref{eg1} and \ref{eg2}. These two theorems, in fact, correspond to the results described in \eqref{etadec} and \eqref{rough} above.
Part II of Section \ref{egsec} contains a typical result -- Theorem \ref{typthm} -- for solutions of equation \eqref{yGeq}. It corresponds to \eqref{Gdec} and \eqref{rough} above.
Moreover, Example \ref{bioapp} shows an application of our theory to a simplistic model in population dynamics. 
Our approach to the analysis of blow-up solutions is compared with the known theory for parabolic equations \cite{GalakBook1995,BeiHu2011} in Remark \ref{comparelit}.

\begin{xnotation}
Throughout the paper, $n\in\N=\{1,2,3,\ldots\}$ is the fixed spatial dimension, i.e., the dimension of the phase space.

For any vector $x\in\R^n$, we denote by $|x|$ its Euclidean norm.

For an $n\times n$ real matrix $A=(a_{ij})_{1\le i,j\le n}$, its Euclidean norm is 
$$\|A\|=\left(\sum_{i=1}^n\sum_{j=1}^n a_{ij}^2\right)^{1/2}.$$

The unit sphere in $\R^n$ is $\mathbb S^{n-1}=\{x\in\R^n:|x|=1\}$.
\end{xnotation}

\section{Existence of finite time blow-ups}\label{blowtime}

We consider the ODE system \eqref{yGeq} in $\R^n$
where the matrix $A$,  functions $H$ and  $G$ will be specified below.

\begin{definition}\label{phom} 
Let $X\ne \{0\}$ and $Y$ be two (real) linear spaces, and $\beta\in\R$ be a given number. 
A function $F:X\setminus\{0\}\to Y$ is positively homogeneous of degree $\beta$ if
\beqs
F(tx)=t^\beta F(x)\text{ for any $x\in X\setminus\{0\}$ and $t>0$.}
\eeqs

Define $\mathcal H_\beta(X,Y)$ to be the set of functions  from $X\setminus\{0\}$ to $Y$ that are positively homogeneous  of degree $\beta$.
\end{definition}

\begin{assumption}\label{assumpAH}
Hereafter,  $A$ is a (real) diagonalizable $n\times n$ matrix with positive eigenvalues,
and $H$ is a function in  $\mathcal H_{\alpha}(\R^n,\R)$ for some $\alpha>0$, and  is positive, continuous on the unit sphere  $\mathbb S^{n-1}$.
\end{assumption}

Thanks to Assumption \ref{assumpAH}, the matrix $A$ has $n$ positive eigenvalues, counting their multiplicities,
\beqs 
\Lambda_1\le \Lambda_2\le \Lambda_3\le \ldots\le \Lambda_n,
\eeqs
and there exists an invertible $n\times n$ (real) matrix $S$ such that
\beq \label{Adiag}
A=S^{-1} A_0 S,
\text{ where } A_0={\rm diag}[\Lambda_1,\Lambda_2,\ldots,\Lambda_n].
\eeq 

Denote the distinct eigenvalues of $A$ by $\lambda_j$'s which are arranged to be (strictly) increasing in $j$, i.e., 
 \beqs
 0<\lambda_1=\Lambda_1<\lambda_2<\ldots<\lambda_{d}=\Lambda_n \text{ for some integer $d\in[1,n]$.}
 \eeqs

The spectrum  of $A$ is $\sigma(A)=\{\Lambda_k:1\le k\le n\}=\{\lambda_j:1\le j\le d\}$.

For $1\le k,\ell\le n$, let $E_{k\ell}$ be the elementary $n\times n$ matrix $(\delta_{ki}\delta_{\ell j})_{1\le i,j\le n}$, where $\delta_{ki}$ and $\delta_{\ell j}$ are the Kronecker delta symbols.
For $\Lambda\in\sigma(A)$, define
 \beqs
\widehat R_\Lambda=\sum_{1\le i\le n,\Lambda_i=\Lambda}E_{ii}\text{ and } R_\Lambda=S^{-1}\widehat R_\Lambda S.
 \eeqs
 
Then one immediately has
\beq \label{Pc} I_n = \sum_{j=1}^{d} R_{\lambda_j},
\quad R_{\lambda_i}R_{\lambda_j}=\delta_{ij}R_{\lambda_j},
\quad  AR_{\lambda_j}=R_{\lambda_j} A=\lambda_j R_{\lambda_j}.
\eeq 

Thanks to \eqref{Pc}, each $R_\Lambda$ is a projection, and $R_\Lambda(\R^n)$ is the eigenspace of $A$ associated with the eigenvalue $\Lambda$.

In the case $A$ is symmetric, the matrix $S$ in \eqref{Adiag} is orthogonal, i.e., $S^{-1}=S^{\rm T}$, 
\beq\label{xAx}
\Lambda_1|x|^2\le x\cdot Ax\le \Lambda_n |x|^2 \text{ for all $x\in\R^n$,}
\eeq
each $R_\Lambda$ is  the orthogonal projection from $\R^n$ to  the eigenspace of $A$ associated with $\Lambda$, and, hence,
\beq \label{Rcontract}
|R_\Lambda x|\le |x| \text{ for  all $x\in\R^n$.}
\eeq 

Regarding the function $H$, we have 
\beqs
0<c_1=\min_{|x|=1} H(x)\le \max_{|x|=1} H(x) =c_2<\infty.
\eeqs 

By writing $ H(x)=|x|^{\alpha} H(x/|x|)$ for any $x\in\R^n\setminus\{0\}$, we derive
\beq\label{Hyal}
c_1 |x|^{\alpha} \le H(x)\le c_2 |x|^{\alpha} \text{ for all }x\in\R^n\setminus \{0\}.
\eeq

In fact, $H$ is continuous and positive on $\R^n\setminus\{0\}$, see (i) and (ii)  of Lemma \ref{Hcts} below.

Let $t_0$ be any given number in $[0,\infty)$. 

 \begin{assumption}\label{assumpF}
   $G(t,x)$ is a continuous function from $[t_0,\infty)\times (\R^n\setminus\{0\})$ to $\R^n$, and 
\beq\label{Grate}
\frac{G(t,x)}{|x|^{1+\alpha}}\to 0 \text{ as $|x|\to\infty$, uniformly in $t\in[t_0,\infty)$.}
\eeq
 \end{assumption}

The followings are standard facts in the theory of ODEs, see, e.g., \cite{JHale78,Hartman1964}.

Consider equation \eqref{yGeq} on the set 
\beq\label{Dset} 
D=\{(t,x):t\ge t_0,x\in \R^n\setminus\{0\}\}\subset \R^{n+1}.
\eeq

For any $y_0\in \R^n\setminus\{0\}$, there exists an interval $[t_0,T_{\max})$, with $0<T_{\max}\le \infty$, and a solution $y\in C^1([t_0,T_{\max}),\R^n\setminus\{0\})$ that satisfies \eqref{yGeq} on $(t_0,T_{\max})$, $y(t_0)=y_0$, and
either
\begin{enumerate}[label=\rnum]
\item $T_{\max}=\infty$, or
\item\label{finT}  $T_{\max}<\infty$, and for any $\varep>0$, $y$ cannot be extended to a function of class $C^1([t_0,T_{\max}+\varep),\R^n\setminus\{0\})$ that satisfies   \eqref{yGeq} on the interval $(t_0,T_{\max}+\varep)$.
\end{enumerate}

It is well-known that, in the case \ref{finT}, it holds, for any compact set $U\subset \R^n\setminus\{0\}$, that
\beq\label{escape}
y(t)\not\in U\text{ when $t\in[t_0, T_{\max})$ is near $T_{\max}$.}
\eeq

Note that such a solution $y(t)$ may not be unique.

\begin{theorem}\label{blowthm}
There exists a number $r_0>0$ such that for any $y_0\in \R^n$ with $|y_0|\ge r_0$ and any solution $y(t)$ on $[t_0,T_{\max})$ as described from \eqref{Dset} to \eqref{escape} above, one has 
\beq\label{limax}
T_{\max}<\infty\text{ and }
 \lim_{t\to T_{\max}^-}|y(t)|=\infty.
\eeq
\end{theorem}
\begin{proof}
Consider $A$ is symmetric first. 
On $(t_0,T_{\max})$, we have
\beq\label{dya1}
\ddt (|y|^{-\alpha})
=-\alpha |y|^{-\alpha-2}y'\cdot y
= -\alpha |y|^{-\alpha-2} H(y)(Ay)\cdot y -\alpha  |y|^{-\alpha-2} G(t,y)\cdot y.
\eeq

Let $a_0= c_1\Lambda_1/2$.
Thanks to \eqref{Grate}, there exists $r_*>0$ such that 
\beq \label{Fcond}
 \frac{|G(t,x)|}{|x|^{1+\alpha}}  \le a_0 \text{ for all  $t \ge t_0$, and all $x \in \R^n$  with $|x|\ge r_*$.}
\eeq 

Set $r_0=4r_*$. For $t>t_0$ sufficiently close to $t_0$, we have $|y(t)|>r_0/2=2r_*$. 
Let $[t_0,T)$ be the maximal interval in $[t_0,T_{\max})$ on which $|y(t)|>r_0/2$.

Suppose $T<T_{\max}$. On the one hand, it must hold that
\beq\label{yT} 
|y(T)|=r_0/2.\eeq 

 On the other hand, we have from \eqref{xAx}, \eqref{Hyal}, \eqref{Fcond}and \eqref{dya1} that, for $t\in(t_0,T)$,
\beq\label{da0}
\ddt (|y|^{-\alpha})
\le -\alpha  c_1\Lambda_1 +\alpha |y|^{-\alpha-1}|G(t,y)|
\le -\alpha  c_1\Lambda_1  + \alpha a_0= -\alpha a_0<0.
\eeq

Thus, $|y(t)|^{-\alpha}\le |y_0|^{-\alpha}$, which implies $|y(t)|\ge |y_0|\ge r_0$ for all $t\in[t_0,T)$. This contradicts \eqref{yT}. Therefore, $T=T_{\max}$. 
Integrating \eqref{da0} gives
\beq\label{ydies}
|y(t)|^{-\alpha}\le |y_0|^{-\alpha}-\alpha a_0( t-t_0)  \text{ for all $t\in[t_0,T_{\max})$.}
\eeq

\medskip
Consider the general matrix $A$ now. Using the equivalence \eqref{Adiag}, we set $\widetilde y(t)=Sy(t)$ and $\widetilde y_0=\widetilde y(t_0)=Sy_0$. 
Then
\beq\label{zeq0}
\ddt \widetilde y=\widetilde H(z)A_0\widetilde y+\widetilde G(t,\widetilde y)\text{ for } t\in(t_0,T_*),
\eeq
where 
\beqs
\widetilde H(z)=H(S^{-1}z),\quad 
\widetilde G(t,z)=SG(t,S^{-1}z) \text{ for $t\in[t_0,T_*)$ and $z\in\R^n\setminus\{0\}$.}
\eeqs

Note that
\beq\label{Syz}
\|S^{-1}\|^{-1}\cdot |x|\le |Sx|\le \|S\|\cdot|x| \text{ for all $x\in\R^n$.}
\eeq

One can verify that $\widetilde H$ belongs to $\mathcal H_{\alpha}(\R^n,\R)$, and, thanks to Lemma \ref{Hcts} below, is positive and continuous on $\R^n\setminus\{0\}$.
Therefore, $\widetilde H$ satisfies the same condition as $H$ in Assumption \ref{assumpAH}. 
Moreover, $\widetilde G(t,z)$ is continuous on $[t_0,\infty)\times (\R^n\setminus\{0\})$.
For $t\in[t_0,T_*)$ and vector $z\ne 0$, we have
\beqs
\frac{|\widetilde G(t,z)|}{|z|^{\alpha+1}}
\le \frac{\|S\| \cdot |G(t,S^{-1}z)|} { \|S^{-1}\|^{-\alpha-1} \cdot  |S^{-1}z|^{\alpha+1}},
\eeqs
which, thanks to \eqref{Grate}, goes to zero as $|z|\to\infty$, uniformly in $t\in[t_0,\infty)$.
Thus, $\widetilde G$ satisfies the same condition as $G$ in Assumption \ref{assumpF}.

Also with $A_0$ being symmetric,  we can apply estimate \eqref{ydies}  to the solution $\widetilde y(t)$ of \eqref{zeq0}. Specifically, there is $r_0>0$ sufficiently large such that when $|y_0|\ge r_0$,  we have $|\widetilde y_0|$ is sufficiently large and
\beqs
|\widetilde y(t)|^{-\alpha}\le |\widetilde y_0|^{-\alpha}-\alpha \widetilde a_0 (t-t_0), \text{ for all $t\in[t_0,T_{\max})$ and some constant $ \widetilde a_0>0$.}
\eeqs

Therefore, 
\beq\label{yagain}
|y(t)|^{-\alpha}\le \|S\|^{-\alpha} |\widetilde y(t)|^{-\alpha} 
\le \|S\|^{-\alpha}(|Sy_0|^{-\alpha}-\alpha \widetilde a_0( t-t_0)) \text{ for all $t\in[t_0,T_{\max})$.}
\eeq

\medskip
If $T_{\max}=\infty$, then \eqref{yagain} implies that $|y(t)|^{-\alpha}<0 $ for $t>t_0+ |Sy_0|^{-\alpha}/(\alpha \widetilde a_0)$, which is an obvious contradiction. Therefore,  $T_{\max}<\infty$. 
As a consequence of \eqref{yagain},  
\beq\label{yunib} 
|y(t)|\ge R_0 \text{ on $[t_0,T_{\max})$, where $R_0=\|S\|\cdot |Sy_0|>0$.}
\eeq

For any $R>R_0$, let $U=\{x\in\R^n:  R_0/2\le |x|\le  R\}$ in \eqref{escape}.
Taking into account \eqref{yunib}, one must have $|y(t)|>R$ when $t\in[t_0, T_{\max})$ is sufficiently close $T_{\max}$.
This proves  the  limit in \eqref{limax}.
\end{proof}

\section{The reference equations}\label{simeg}

Let $a$ and $\alpha$ be arbitrarily  positive numbers and $t_0,T_*\in\R$ be two given numbers with $T_*>t_0\ge 0$. 
Assume $y\in C^1([t_0,T_*),\R^n)$ satisfies 
\beq \label{ynozero}
\text{$y(t)\ne 0$ for all $t\in [t_0,T_*)$, }
\eeq
\beq \label{limbig}
\lim_{t\to T_*^-} |y(t)|= \infty,
\eeq 
and 
\beq\label{basicf}
y'=a|y|^{\alpha} y +f(t) \text{ for } t \in (t_0,T_*),
\eeq
where $f:[t_0,T_*)\to \R^n$ is continuous and satisfies
\beq\label{fE}
|f(t)|\le |y(t)|^{1+\alpha}E_0(t) \text{ for all $t\in[t_0,T_*)$} 
\eeq
with some 
\beq \label{E0ic}
\text{Lebesgue integrable function } E_0:[t_0,T_*)\to [0,\infty).
\eeq 

We start with a derivation of  both lower and upper bounds  for $|y(t)|$.

\begin{lemma}\label{simlem}
Suppose
\beq\label{eavg}
e_*(t)\eqdef \frac{1}{T_*-t}\int_t^{T_*}E_0(\tau)\d\tau \to 0\text{ as $t\to T_*^-$.}
\eeq
Then there exist $C_1,C_2>0$ such that
\beq\label{ysimest}
C_1(T_*-t)^{-1/\alpha}
\le  |y(t)|\le 
 C_2 (T_*-t)^{-1/\alpha} \text{ for all $t\in[t_0,T_*)$.}
\eeq
\end{lemma}
\begin{proof}
For $t\in( t_0,T_*)$, taking into account \eqref{ynozero}, we have 
\beq\label{dy3}
\ddt(|y|^{-\alpha})=-\alpha a -\alpha |y|^{-\alpha-2}f(t)\cdot y.
\eeq

Integrating equation \eqref{dy3} from $t$ to $T_*$ and using the limit \eqref{limbig} (or, more precisely, integrating \eqref{dy3} from $t$ to $t'\in(t,T_*)$ and letting $t'\to T_*^-$) give
\beq\label{ynorm}
|y(t)|^{-\alpha}=\alpha a (T_*-t)+g(t), 
\text{ where } g(t)=\alpha \int_t^{T_*} |y(\tau)|^{-\alpha-2}f(\tau)\cdot y(\tau)\d\tau.
\eeq

By the Cauchy--Schwarz inequality and \eqref{fE},  one has
$|y(\tau)|^{-\alpha-2}|f(\tau)\cdot y(\tau)|\le E_0(\tau)$. Because $E_0\in L^1(0,T_*)$, the $g(t)$ in \eqref{ynorm}, indeed, is a continuous function from $[t_0,T_*)$ to $\R$.
Moreover,
\beq\label{gest}
|g(t)|\le \alpha \int_t^{T_*} E_0(\tau) \d\tau
=\alpha(T_*-t)e_*(t) . 
\eeq

By \eqref{eavg}, there is $T\in [t_0,T_*)$ such that $e_*(t)\le a/2$ for all $t\in [T,T_*)$. Combining this with \eqref{ynorm} and \eqref{gest} gives, for all $t\in [T,T_*)$,
\beqs
(\alpha a/2) (T_*-t)\le |y(t)|^{-\alpha}\le (3\alpha a/2) (T_*-t).
\eeqs
Thus,
\beq\label{yratesim}
\left (\frac2{3\alpha a}\right)^{1/\alpha}(T_*-t)^{-1/\alpha}
\le  |y(t)|\le 
 \left(\frac2{\alpha a}\right)^{1/\alpha} (T_*-t)^{-1/\alpha} \text{ for all $t\in[T,T_*)$.}
\eeq

Combining \eqref{yratesim} with the fact
 \beq\label{yTT}
0<\min_{t\in[t_0,T]} (T_*-t)^{1/\alpha}|y(t)| \le \max_{t\in[t_0,T]} (T_*-t)^{1/\alpha}|y(t)|<\infty,
\eeq
we obtain \eqref{ysimest}.
\end{proof}

Note that if $E_0(t)\to 0$ as $t\to T_*^-$, then condition \eqref{eavg} is met.
The main result of this section is the next theorem.

\begin{theorem}\label{simthm}
Assume 
\beq\label{condE}
Z_1\eqdef \int_{t_0}^{T_*}\frac{E_0(\tau)}{T_*-\tau}\d\tau <\infty.
\eeq

Then $ (T_*-t)^{1/\alpha} y(t)$ converges to a nonzero vector $\xi_*$  as $t\to T_*^-$.
More specifically,
\beq\label{yxi}
| (T_*-t)^{1/\alpha} y(t)-\xi_* |=\bigo(E_1(t))\text{ as $t\to T_*^-$,}
\eeq
where
\beq\label{E1def}
E_1(t)=\int_t^{T_*}\frac{E_0(\tau)}{T_*-\tau}\d\tau\text{  for $t\in[t_0,T_*)$,}
\eeq
and $\xi_*$ satisfies
\beq\label{axi1}
 \alpha a |\xi_*|^\alpha=1.
\eeq
\end{theorem}
\begin{proof}
Note that $E_1(t)$ is nonnegative, decreasing on $[t_0,T_*)$, and 
\beq\label{Elim}
\lim_{t\to T_*^-} E_1(t)=0.
\eeq
Hence, it suffices to prove \eqref{yxi} and \eqref{axi1}.
It is clear  that 
\beq \label{ese1}
e_*(t)\le E_1(t) \text{ for all $t\in[t_0,T)$.}
\eeq 

Because of \eqref{ese1} and \eqref{Elim}, $e_*(t)$ satisfies \eqref{eavg}, and the estimates in \eqref{ysimest} hold true. We use the calculations in the proof of Lemma \ref{simlem}.
By \eqref{gest} and \eqref{ese1},
\beq\label{gest2}
|g(t)|\le \alpha  (T_*-t)E_1(t) . 
\eeq

Considering equation \eqref{basicf} as a linear equation of $y$,  we solve for $y(t)$ explicitly  by the variation of constants formula. It results in
\begin{align*}
y(t)&=e^{J(t)}\left(y_0+\int_{t_0}^t e^{-J(\tau)}f(\tau)\d\tau\right) \text{ for $t\in [t_0,T_*)$,}
\end{align*}
where
\beq\label{Iform}
J(t)=a\int_{t_0}^t |y(\tau)|^{\alpha}\d\tau.
\eeq

As a consequence of \eqref{ynorm},  one has $\alpha a (T_*-t)+g(t)>0$ for all $t\in[t_0,T_*)$.
Combining \eqref{ynorm} with \eqref{Iform}, we split $J(t)$ into 
\begin{align*}
J(t)=\int_{t_0}^t \frac{a}{a\alpha (T_*-\tau)+  g(\tau)}\d\tau=J_1(t)+J_2(t),
\end{align*}
where
\beqs
J_1(t)=\int_{t_0}^t \frac{1}{\alpha (T_*-\tau)}\d\tau
\text{ and }
J_2(t)=\int_{t_0}^t h_0(\tau) \d\tau,
\eeqs 
with
\beqs 
h_0(\tau)=\frac{ -g(\tau)}{\alpha (T_*-\tau)(a\alpha (T_*-\tau)+ g(\tau))}.
\eeqs

It is obvious that
$e^{J_1(t)}=(T_*-t_0)^{1/\alpha}(T_*-t)^{-1/\alpha}$.
Therefore,
\beq\label{yJform}
y(t)=\frac{(T_*-t_0)^{1/\alpha}}{(T_*-t)^{1/\alpha}} e^{J_2(t)}\left(y_0+(T_*-t_0)^{-1/\alpha}\int_{t_0}^t e^{-J_2(\tau)}(T_*-\tau)^{1/\alpha}f(\tau)\d\tau\right).
\eeq

Using  \eqref{gest2} to estimate $|g(\tau)|$ in the denominator of $h_0(t)$, and then using \eqref{gest} to estimate the remaining $|g(t)|$, we find that, as $\tau\to T_*^-$,
\beq\label{htau}
|h_0(\tau)|=\bigo(|g(\tau)|(T_*-\tau)^{-2})=\bigo(e_*(\tau)(T_*-\tau)^{-1}). 
\eeq

By Fubini's theorem, we calculate 
\beq\label{e2E}
 \int_t^{T_*} e_*(\tau)(T_*-\tau)^{-1} \d\tau
 = \int_t^{T_*}  E_0(s)\left(\frac1{T_*-s}-\frac1{T_*-t}\right)\d s
 \le E_1(t)<\infty.
\eeq

Estimates \eqref{htau} and \eqref{e2E} imply 
\beq\label{J2lim} 
\lim_{t\to T_*^-}J_2(t)=\int_{t_0}^{T_*} h_0(\tau)\d\tau=J_*\in\R.
\eeq 

Combining  \eqref{J2lim}, \eqref{fE} and \eqref{ysimest}, we have 
\beq\label{eJf}
e^{J_2(t)} |f(t) |(T_*-t)^{1/\alpha}
=\bigo((T_*-t)^{-1}E_0(t))  \text{ as $t\to T_*^-$.}
\eeq

By the assumption \eqref{condE},
\beqs
\lim_{t\to T_*^-} \int_{t_0}^t e^{J_2(\tau)} (T_*-\tau)^{1/\alpha}f(\tau)\d\tau
=\int_{t_0}^{T_*} e^{J_2(\tau)}(T_*-\tau)^{1/\alpha}f(\tau)\d\tau
=\eta_*\in \R^n.
\eeqs

Denote 
\beq\label{JJh1}
h_1(t)=J_*-J_2(t)=\int_t^{T_*}h_0(\tau)\d\tau\in \R
\eeq
and 
\beq\label{ifeeta}
\eta(t)=\eta_*-\int_{t_0}^t e^{J_2(\tau)}(T_*-\tau)^{1/\alpha}f(\tau)\d\tau
=\int_t^{T_*} e^{J_2(\tau)}(T_*-\tau)^{1/\alpha}f(\tau)\d\tau\in \R^n.
\eeq

Combining \eqref{yJform}, \eqref{JJh1} and \eqref{ifeeta} gives, for $t\in[t_0, T_*)$,
\begin{align*}
(T_*-t)^{1/\alpha}y(t)
&=e^{J_2(t)} \left((T_*-t_0)^{1/\alpha} y_0+\eta_*-\eta(t)\right)\\
&=e^{J_*-h_1(t)} \left((T_*-t_0)^{1/\alpha} y_0+\eta_*\right)-e^{J_2(t)}\eta(t)
=e^{-h_1(t)} \xi_*- e^{J_2(t)}\eta(t),
\end{align*}
where $\xi_*=e^{J_*}( (T_*-t_0)^{1/\alpha}y_0+\eta_*)\in\R^n$.
It follows that
\begin{align*}
\left|(T_*-t)^{1/\alpha} y(t)-\xi_* \right|
=\bigo\left(|e^{-h_1(t)}-1|+|\eta(t)|\right),
\end{align*}

By estimates \eqref{htau} and \eqref{e2E}, we have, for $t$ close to $T_*$, 
\beqs
|h_1(t)|\le C\int_t^{T_*} e_*(\tau)(T_*-\tau)^{-1} \d\tau\le C E_1(t)
\text{ for some $C>0$.}
\eeqs

It follows \eqref{eJf} that, for $t$ close to $T_*$, 
\beqs
|\eta(t)| \le C' \int_t^{T_*} (T_*-\tau)^{-1}E_0(\tau) \d\tau=C' E_1(t)
\text{ for some $C'>0$.}
\eeqs

Therefore,
\beq \label{yy1}
\left|(T_*-t)^{1/\alpha} y(t)-\xi_* \right|=\bigo\left(|h_1(t)|+|\eta(t)|\right)
=\bigo(E_1(t)),
\eeq
which proves \eqref{yxi}. 

We  prove \eqref{axi1} next.
  On the one hand, one has from \eqref{ynorm} that
\beqs
(T_*-t)^{1/\alpha}|y(t)|= \left(a\alpha+\frac{g(t)}{ T_*-t}\right)^{-1/\alpha}.
\eeqs
  
By properties \eqref{gest2} and  \eqref{Elim},  we have 
$\lim_{t\to T_*^-} g(t)/( T_*-t)=0$.
Thus,
\beq\label{yxinorm2}
\lim_{t\to T_*^-} (T_*-t)^{1/\alpha}|y(t)| =(a\alpha)^{-1/\alpha}.
\eeq

On the other hand, we have from  \eqref{yy1} that
  \beq\label{yxinorm1}
\lim_{t\to T_*^-}(T_*-t)^{1/\alpha} |y(t)|= |\xi_*|.
\eeq

By the limits \eqref{yxinorm1} and \eqref{yxinorm2}, we obtain $|\xi_*|=(a\alpha)^{-1/\alpha}$, which yields \eqref{axi1}.
\end{proof}
 
 We note that condition \eqref{condE}, in fact,  already implies \eqref{E0ic}.

\section{The rate of blow-ups}\label{ratesec}
 
Let $t_0,T_*\in\R$ be two given numbers with $T_*>t_0\ge 0$.
Assume $y\in C^1([t_0,T_*),\R^n)$ has properties \eqref{ynozero}, \eqref{limbig}, and
satisfies equation \eqref{mainode} in the interval $(t_0,T_*)$,
where the function $f:[t_0,T_*)\to \R^n$ is continuous and satisfies \eqref{fE}, \eqref{E0ic}.

We establish a counterpart of Lemma \ref{simlem}.

\begin{lemma}\label{estthm}
Assume \eqref{eavg}. Then there are positive constants  $C_1$ and $C_2$ such that 
 \beq\label{yplus}
 C_1(T_*-t)^{-1/\alpha}
\le  |y(t)|\le 
 C_2(T_*-t)^{-1/\alpha} \text{ for all $t\in[t_0,T_*)$.}
\eeq
 \end{lemma}
\begin{proof}
We treat the cases of symmetric $A$  and non-symmetric $A$ separately.

\medskip
\noindent\textit{Part 1.} Consider $A$ is symmetric. 
The proof in this case is similar to that of Lemma \ref{simlem}.
For $t\in(t_0,T_*)$, we calculate
\beqs
\ddt (|y|^{-\alpha})
= -\alpha |y|^{-\alpha-2} H(y)(Ay)\cdot y -\alpha  |y|^{-\alpha-2} f(t)\cdot y.
\eeqs

Utilizing  \eqref{xAx}, \eqref{Hyal}  and \eqref{fE}, one has 
\beqs
-\alpha c_2 \Lambda_n   -\alpha |y|^{-\alpha-2} f(t)\cdot y\le \ddt  (|y|^{-\alpha})\le  -\alpha  c_1 \Lambda_1  -\alpha |y|^{-\alpha-2} f(t)\cdot y.
\eeqs

Integrating from $t$ to $T_*$ gives
\beqs
\alpha  c_1 \Lambda_1 (T_*-t)  +g(t)\le |y(t)|^{-\alpha}\le \alpha c_2 \Lambda_n(T_*-t) +g(t),
\eeqs
where $g(t)$ is defined in \eqref{ynorm}.

Let $\varep_0\in (0,c_1\Lambda_1)$. By \eqref{gest} and \eqref{eavg},  there is $T\in[t_0,T_*)$ such  that $|g(t)|\le \alpha \varep_0 (T_*-t)$
on $[T,T_*)$.
Hence, 
\beqs
\alpha (c_a\Lambda_1-\varep_0) (T_*-t)\le  |y(t)|^{-\alpha} \le \alpha (c_2 +\varep_0) (T_*-t) \text{ for all $t\in [T,T_*)$,}
 \eeqs
which yields
 \beq\label{ytclose}
( \alpha(c_2\Lambda_n+\varep_0)(T_*-t))^{-1/\alpha}
\le  |y(t)|\le 
(\alpha(c_1\Lambda_1-\varep_0)(T_*-t))^{-1/\alpha} \text{ for all $t\in[T,T_*)$.}
\eeq

Fixing $\varep_0=c_1\Lambda_1/2$ in \eqref{ytclose} and combining it with the fact \eqref{yTT}, we obtain the desired estimates in \eqref{yplus}. 

\medskip
\noindent\textit{Part 2.}
Consider $A$ is non-symmetric.
Let $A=S^{-1}A_0 S$ as in \eqref{Adiag}.

Set $\widetilde y(t)=Sy(t)$ for $t\in[t_0,T_*)$.  
Then $\widetilde y$ belongs to $C^1([t_0,T_*),\R^n)$, 
$\widetilde y(t)\ne 0$ for all $t\in[t_0,T_*)$, $\lim_{t\to T_*^-}|\widetilde y(t)|=\infty$, and 
\beq\label{zeq}
\ddt \widetilde y=\widetilde H(\widetilde y)A_0\widetilde y+\widetilde{f}(t)\text{ for } t\in(t_0,T_*),
\eeq
where 
\beq\label{HRz}
\widetilde H(z)=H(S^{-1}z) \text{ for $z\in\R^n\setminus\{0\}$, and }
\widetilde {f}(t)=Sf(t) \text{ for $t\in[t_0,T_*)$.}
\eeq

Same as in the proof of Theorem \ref{blowthm}, $\widetilde H$ belongs to $\mathcal H_{\alpha}(\R^n,\R)$,  is positive and continuous on $\R^n\setminus\{0\}$.

Clearly, the function $\widetilde f$ is continuous on $[t_0,T_*)$. Thanks to \eqref{fE} and \eqref{Syz}, it satisfies, for $t\in[t_0,T_*)$, 
\beqs
|\widetilde f(t)|\le \|S\|\cdot |f(t)|\le \|S\| \cdot  |y(t)|^{1+\alpha}E_0(t)\le |\widetilde y(t)|^{1+\alpha}\widetilde E_0(t),
\eeqs
where
\beq\label{tilEzero}
\widetilde E_0(t)=\|S\| \cdot \|S^{-1}\|^{\alpha+1}E_0(t).
\eeq

Combining \eqref{tilEzero} with  \eqref{E0ic} and \eqref{eavg} yields  
\beqs
\widetilde E_0\in L^1(t_0,T_*)\text{ and }\lim_{t\to T_*^-} \frac{1}{T_*-t}\int_t^{T_*}\widetilde E_0(\tau)\d\tau=0.
\eeqs

Therefore, with $A_0$ being symmetric, we can apply the result in Part 1 to  the solution $\widetilde y(t)$ and equation \eqref{zeq}. Then there exist two positive constants $\widetilde C_1$ and $\widetilde C_2$ such that 
\beq\label{zpower}
\widetilde C_1(T_*-t)^{-1/\alpha}
\le  |\widetilde y(t)|\le 
 \widetilde C_2(T_*-t)^{-1/\alpha}
\text { for all $t\in [t_0,T_*)$.}
 \eeq

Combining \eqref{zpower} with the relations in \eqref{Syz}, we obtain the estimates in \eqref{yplus} for $y(t)$. 
\end{proof}

By \eqref{Hyal} and \eqref{yplus}, we immediately have in Lemma \ref{estthm},  for all $t\in[t_0,T_*)$, that
\beq\label{Hyt}
C_3 (T_*-t)^{-1}\le H(y(t)) \le C_4 (T_*-t)^{-1}, \text{ where $C_3=c_1C_1^{\alpha}$ and $C_4=c_2C_2^{\alpha}$.} 
\eeq

From Lemma \ref{estthm}, we can derive the lower and upper bounds for blow-up solutions of equation \eqref{yGeq} as follows.  

\begin{corollary}\label{boundcor}
Let $G(t,x)$ be as in Assumption \ref{assumpF}.
Suppose $y\in C^1([t_0,T_*),\R^n\setminus\{0\})$  satisfies equation \eqref{yGeq} on $(t_0,T_*)$ and has property \eqref{limbig}. Then the estimates in \eqref{yplus} hold true. 
\end{corollary}
\begin{proof}
Let  $f(t)=G(t,y(t))$ for $t\in[t_0,T_*)$. We consider $y(t)$ as a solution of equation \eqref{mainode}  with $E_0(t)$ in \eqref{fE} being 
$$E_0(t)=\frac{|f(t)|}{|y(t)|^{\alpha+1}}=\frac{|G(t,y(t))|}{|y(t)|^{\alpha+1}}.$$

Then $E_0(t)$ is continuous on $[t_0,T_*)$ and, thanks to \eqref{Grate} and \eqref{limbig}, $E_0(t)\to 0$ as $t\to T_*^-$. 
Thus, $f$ and $E_0$ satisfy \eqref{fE}, \eqref{E0ic} and \eqref{eavg}. Applying Lemma \ref{estthm},   we obtain the estimates in \eqref{yplus}.
\end{proof}

\section{The symmetric matrix case}\label{symcase}

In this section, we assume the matrix $A$ is symmetric.
Let $y(t)$ be a solution of equation \eqref{mainode}  as in the beginning of Section \ref{ratesec}. 
For $t\in [t_0,T_*)$, define
\beq\label{quot}
\lambda(t)=\frac{y(t)\cdot Ay(t)}{|y(t)|^2} \text{ and }   
 v(t)=\frac{y(t)}{|y(t)|} .
\eeq

Then $\lambda\in C^1([t_0,T_*),\R)$ and   $v\in C^1([t_0,T_*),\R^n)$. Moreover, one has, 
$|v(t)|=1$ and, thanks to \eqref{xAx},  
 \beq\label{lest}
 \Lambda_1 \le  \lambda(t)\le \Lambda_n\le \|A\| \text{ for all $t\in[ t_0,T_*)$.}
 \eeq 
 
 We will find asymptotic behavior of $\lambda(t)$, then of $v(t)$, and finally of $y(t)$.
 
 For the rest of this section, we assume \eqref{condE} and define  $E_1(t)$ by \eqref{E1def}. Then Lemma \ref{estthm} applies and the estimates in \eqref{yplus} hold true.
 
\begin{proposition}\label{lem1}
One has 
$$\lim_{t\to T_*^-} \lambda(t)=\Lambda\in \sigma(A).$$
\end{proposition}
\begin{proof}
For $t\in(t_0,T_*)$, we have
\beqs
\lambda'(t)
=\frac2{|y|^2} y' \cdot (Ay-\lambda y)
=\frac2{|y|^2}\Big [H(y)(Ay-\lambda y) + \lambda H(y)y+f(t)\Big]\cdot (Ay-\lambda y).
\eeqs

Because $ y\cdot(Ay-\lambda y)=0$, it follows that 
\beq\label{lh}
\lambda'(t)
=2H(y)|Av-\lambda v|^2 +h(t),
\eeq
where
\beqs
h(t)=\frac2{|y(t)|^2}f(t)\cdot(Ay(t)-\lambda(t) y(t)).
\eeqs

Using \eqref{fE}, \eqref{lest}   and, then,  \eqref{yplus}, we estimate  
\beq\label{oh}
|h(t)|\le 4\|A\|\cdot |y(t)|^{\alpha} E_0(t)\le C_5 (T_*-t)^{-1} E_0(t)
\text{ for all $t\in [t_0,T_*)$,}
\eeq
 where $C_5=4\|A\| C_2^{\alpha}$. 

For $t,t'\in[t_0,T_*)$ with $t'>t$, integrating equation \eqref{lh} from $t$ to $t'$ gives
\beq\label{llH}
\lambda(t')-\lambda(t)=2\int_t^{t'}H(y(\tau))|Av(\tau)-\lambda(\tau) v(\tau)|^2\d\tau +\int_t^{t'}h(\tau)\d\tau.
\eeq

Thanks to \eqref{oh}, the last integral can be estimated as
\beq\label{hitt}
\left|\int_t^{t'}h(\tau)\d\tau\right|\le C_5E_1(t).
\eeq

By taking the limit inferior of \eqref{llH}, as $t'\to T_*^-$, we derive
\beq\label{llH2}
\liminf_{t'\to T_*^-}\lambda(t')\ge \lambda(t)-C_5 E_1(t).
\eeq

Then taking the limit superior of \eqref{llH2}, as $t\to T_*^-$, yields
\beqs
\liminf_{t'\to T_*^-}\lambda(t')\ge \limsup_{t\to T_*^-}  \lambda(t).
\eeqs

This and \eqref{lest}  imply
\beq\label{limlam} 
\lim_{t\to T_*^-}\lambda(t)=\Lambda\in [\Lambda_1,\Lambda_n].
\eeq

Using properties \eqref{hitt} and  \eqref{limlam} in \eqref{llH} and by taking $t=t_0$ and $t'\to T_*^-$, we obtain
\beq\label{Hfin}
\int_{t_0}^{T_*} H(y(\tau))|Av(\tau)-\lambda(\tau) v(\tau)|^2\d\tau<\infty.
\eeq

We claim that 
\beq \label{claim} 
\forall\varep\in(0,T_*-t_0),\exists t\in [T_*-\varep,T_*):|Av(t)-\lambda(t)v(t)|<\varep.
\eeq 

Indeed, suppose the claim \eqref{claim} is not true, then 
\beq\label{anticlaim} 
\exists\varep_0\in(0,T_*-t_0),\forall t\in [T_*-\varep_0,T_*): |Av(t)-\lambda(t)v(t)|\ge \varep_0.
\eeq

Combining \eqref{anticlaim} with property \eqref{Hyt}, we have
\beqs
\int_{T_*-\varep_0}^{T_*} H(y(\tau))|Av(\tau)-\lambda(\tau) v(\tau)|^2\d\tau
\ge \int_{T_*-\varep_0}^{T_*} C_3(T_*-\tau)^{-1} \varep_0^2 \d\tau=\infty,
\eeqs
which contradicts \eqref{Hfin}. Hence, the  claim \eqref{claim} is true. 

Thanks to \eqref{claim}, there exists a sequence $(t_j)_{j=1}^\infty\subset [t_0,T_*)$ such that  
\beq \label{Avj}
\lim_{j\to\infty}t_j= T_* 
\text{ and }
\lim_{j\to\infty}|Av(t_j)-\lambda(t_j)v(t_j)|= 0.
\eeq 

The first equation in \eqref{Avj} and \eqref{limlam} imply $\lambda(t_j)\to \Lambda$ as $j\to\infty$. 
Because $v(t_j)\in\mathbb S^{n-1}$ for all $j$, we can extract a subsequence $(v(t_{j_k}))_{k=1}^\infty$ such that $v(t_{j_k})\to \bar v\in\mathbb S^{n-1}$ as $k\to\infty$. 
Combining these limits with the second equation in \eqref{Avj} written with $j=j_k$ and $k\to\infty$ yields
$A\bar v=\Lambda\bar v$.
Therefore, $\Lambda$ is an eigenvalue of $A$.
\end{proof}

\textit{In the remainder  of this section, $\Lambda$ is the eigenvalue in Proposition \ref{lem1}.}

\medskip
For the asymptotic behavior of $v(t)$ as $t\to T_*^-$, we investigate that of $(I_n-R_{\Lambda})v(t)$ first, and then of $R_{\Lambda}v(t)$.

If $\Lambda_1=\Lambda_n$, then $\sigma(A)=\{\Lambda\}$,  thus,
\beq\label{ILv}
 I_n=R_{\Lambda}
\text{ and } (I_n-R_{\Lambda})v(t)=0.
\eeq

\begin{assumption}\label{E2assum}
In the case  $\Lambda_1<\Lambda_n$, we assume further that
\beq\label{condE2}
Z_2 \eqdef  \int_{t_0}^{T_*} (T_*-\tau)^{-1}E_0^2(\tau)\d\tau<\infty,
\eeq
and define
\begin{align}  \label{mudef}
\mu&=\min \{ |\lambda_j-\Lambda| : 1\le j\le d, \lambda_j\ne  \Lambda\},\\
\label{E2def}
E_2(t)&=\left(\int_{t}^{T_*} (T_*-\tau)^{-1}E_0^2(\tau)\d\tau\right)^{1/2} \text{ for $t\in[t_0,T_*)$.}
\end{align}
\end{assumption}

In Assumption \ref{E2assum}, one has $\mu>0$ and, thanks to \eqref{condE2}, 
\beq\label{limE2}
\lim_{t\to T_*^-}E_2(t)=0.
\eeq

\begin{proposition}\label{lem2}
One has
\beq\label{remvstar}
\lim_{t\to T_*^-} (I_n-R_{\Lambda})v(t)=0.
 \eeq

More specifically, by setting
\beq\label{Ebarate}
V_1(t)= |(I_n-R_\Lambda)v(t)|,
\eeq
 the following estimates hold.

If $\Lambda_1=\Lambda_n$, then 
\beq \label{rvstar0}
V_1(t)=0 \text{  for all $t\in[t_0,T_*)$.}
\eeq

If $\Lambda_1=\Lambda<\Lambda_n$,  then
\beq\label{rvstar1}
V_1(t)=\bigo(E_2(t))  \text{  as $t\to T_*^-$.}
 \eeq
 
 If  $\Lambda_1<\Lambda\le \Lambda_n$, then there is $\theta_0>0$ such that, for any $\varep\in(0,1)$, 
 \beq
 \label{rvstar2}
V_1(t)=\bigo((T_*-t)^{(1-\varep)\theta_0}+E_2(T_*-(T_*-t)^\varep))
\text{  as $t\to T_*^-$.}
\eeq
 \end{proposition}
\begin{proof}
When $\Lambda_1=\Lambda_n$, \eqref{rvstar0} follows \eqref{ILv}.
Consider $\Lambda_1<\Lambda_n$, which gives $\sigma(A)\ne \{\Lambda\}$.
We calculate
\beqs
 v' =\frac1{|y|}y' -\frac1{|y|^3}(y' \cdot y)y
 =\frac{H(y)}{|y|}Ay+\frac1{|y|}f(t)
-\frac{H(y)(Ay)\cdot y}{|y|^3}y
 -\frac{f(t)\cdot y}{|y|^3}y.
\eeqs

Define the function $g:[t_0,T_*)\to\R^n$  by
\beqs
g(t)=\frac{f(t)}{|y(t)|}-\left(\frac{f(t)}{|y(t)|}\cdot v(t)\right)v(t).
\eeqs

Then we have
\beq \label{dv}
 v'=H(y)(Av-\lambda v)+g(t)\text{ for all $t\in (t_0,T_*)$.}
\eeq 

Using \eqref{fE}, we estimate
\beq\label{gMy}
|g(t)|\le 2 |f(t)|/|y(t)| \le 2 |y(t)|^{\alpha}E_0(t)\text{ for all $t\in[t_0,T_*)$.}
\eeq

Let $\lambda_j\in\sigma(A)\setminus\{\Lambda\}$.
Applying $R_{\lambda_j}$ to equation \eqref{dv} and taking the dot product with $R_{\lambda_j}v$ yield
\beq\label{Rnorm}
\frac12 \ddt |R_{\lambda_j}v|^2
 =H(y)(\lambda_j-\lambda)| R_{\lambda_j}v|^2 +R_{\lambda_j} g(t)\cdot R_{\lambda_j}v.
\eeq

Applying  Cauchy--Schwarz's inequality, inequality \eqref{Rcontract} to $|R_{\lambda_j}g(t)|$,   and then Cauchy's inequality, we have
\beqs
|R_{\lambda_j}g(t)\cdot R_{\lambda_j}v|
\le |g(t)||R_{\lambda_j}v|
\le \frac\mu4 H(y)| R_{\lambda_j}v|^2 + \frac{|g(t)|^2}{\mu H(y)}.
\eeqs

Applying the first inequality of \eqref{Hyal} to estimate the last $H(y)$ and using \eqref{gMy} give 
\beqs
|R_{\lambda_j}g(t)\cdot R_{\lambda_j}v|
\le \frac\mu4 H(y)| R_{\lambda_j}v|^2 + \frac{4}{\mu c_1}|y|^\alpha E_0^2(t).
\eeqs

Utilizing the upper bound of $|y(t)|$ in \eqref{yplus}, we obtain, for $t\in [t_0,T_*)$,
\beq\label{abso}
|R_{\lambda_j}g(t)\cdot R_{\lambda_j}v|
\le \frac\mu4 H(y)| R_{\lambda_j}v|^2 +\frac{C_6}2 (T_*-t)^{-1}E_0^2(t),
\text{ where
$C_6=\frac{8C_2^\alpha}{\mu c_1}$.}
\eeq

Below, $T\in (t_0,T_*)$ is fixed and can be taken sufficiently close to $T_*$ such that
\beq\label{Tlam}
|\lambda(t)-\Lambda|\le \frac\mu4 \text{ for all $t\in[T,T_*)$.}
\eeq

\medskip
\noindent\underline{\textit{Case $\lambda_j>\Lambda$.}} In this case, combining \eqref{Rnorm} and \eqref{abso} yields, for $t\in (t_0,T_*)$,
\begin{align*}
\frac12 \ddt |R_{\lambda_j}v|^2
 &\ge (\lambda_j-\lambda-\frac\mu4) H(y)| R_{\lambda_j}v|^2 -\frac{C_6}2 (T_*-t)^{-1}E_0^2(t) .
\end{align*}

By definition \eqref{mudef} of $\mu$ and the choice \eqref{Tlam}, one has,  for all $t\in [T,T_*)$,
\beq\label{lamu} 
\lambda_j-\lambda(t)-\frac\mu4= (\lambda_j-\Lambda)+(\Lambda-\lambda(t))-\frac\mu4
\ge \mu-\frac\mu4-\frac\mu4=\frac{\mu}2.
\eeq

Hence,
\begin{align*}
\ddt |R_{\lambda_j}v|^2
 &\ge \mu H(y)| R_{\lambda_j}v|^2 -C_6(T_*-t)^{-1}E_0^2(t).
\end{align*}

Then, for any $t,\bar t\in  [T,T_*)$ with $t>\bar t$, one has
\beq\label{eHR}
e^{-\mu\int_{\bar t}^t H(y(\tau))\d\tau}  |R_{\lambda_j}v(t)|^2-  |R_{\lambda_j}v(\bar t)|^2
 \ge -C_6 \int_{\bar t}^t e^{-\mu\int_{\bar t}^\tau H(y(s))\d s}(T_*-\tau)^{-1}E_0^2(\tau)\d\tau.
\eeq

For the first term on the left hand side of \eqref{eHR}, note from the first inequality in \eqref{Hyt}  that $\int_{\bar t}^{T_*} H(y(\tau))\d\tau=\infty$, and from \eqref{Rcontract} that $ |R_{\lambda_j}v(t)|\le |v(t)|=1$. Then letting $t\to T_*^-$ in \eqref{eHR} yields
\begin{align*}
|R_{\lambda_j}v(\bar t)|^2 
 &\le C_6 \int_{\bar t}^{T_*} e^{-\mu\int_{\bar t}^\tau H(y(s))\d s}(T_*-\tau)^{-1}E_0^2(\tau)\d\tau\\
& \le  C_6 \int_{\bar t}^{T_*} (T_*-\tau)^{-1}E_0^2(\tau)\d\tau
= C_6 E_2^2(\bar t).
\end{align*}

Therefore, we obtain
\beq\label{Rv2} 
|R_{\lambda_j}v(\bar t)| =\bigo(E_2(\bar t))\text{ as $\bar t\to T_*^-$.}
\eeq 

\medskip
\noindent\underline{\textit{Case $\lambda_j<\Lambda$.}} Using \eqref{abso} to estimate  the last term in \eqref{Rnorm} again, we have
\begin{align*}
\frac12 \ddt |R_{\lambda_j}v|^2
 &\le (\lambda_j-\lambda+\frac\mu4) H(y)| R_{\lambda_j}v|^2 +\frac{C_6}2 (T_*-t)^{-1}E_0^2(t).
\end{align*}

Similar to \eqref{lamu}, one has, for  $t\in [T,T_*)$, 
$$\lambda_j -\lambda(t)+\frac\mu4=(\lambda_j-\Lambda)+ (  \Lambda-\lambda(t))+\frac\mu4
\le -\mu+\frac\mu4 +\frac{\mu}4=-\frac{\mu}2.$$

Thus, for $t\in [T,T_*)$,
\beq\label{dRlarge}
 \ddt |R_{\lambda_j}v|^2
 \le - \mu H(y)| R_{\lambda_j}v|^2 +C_6(T_*-t)^{-1}E_0^2(t).
\eeq

Let $t$ and $\bar t$ be any numbers in $[T,T_*)$ with  $t>\bar t$.
 It follows \eqref{dRlarge} that 
\beqs
  |R_{\lambda_j}v(t)|^2
 \le e^{-\mu\int_{\bar t}^t H(y(\tau))\d\tau} | R_{\lambda_j}v(\bar t)|^2 + C_6\int_{\bar t}^t e^{-\mu\int_{\tau}^t H(y(s))\d s}(T_*-\tau)^{-1}E_0^2(\tau)\d\tau.
 \eeqs

We simply estimate 
\beq\label{simR}
  |R_{\lambda_j}v(t)|^2
\le e^{-\mu\int_{\bar t}^t H(y(\tau))\d\tau}  + C_6E_2^2(\bar t).
 \eeq 

Utilizing  \eqref{Hyt} to find a lower bound of $H(y(\tau))$ in \eqref{simR}, we find
\beq\label{RE3}
  |R_{\lambda_j}v(t)|^2
\le  \frac{(T_*-t)^{C_3\mu}}{(T_*-\bar t)^{C_3\mu}}  
  + C_6 E_2^2(\bar t).
\eeq

Letting $t\to T_*^-$, we obtain
\beqs
\limsup_{t\to T_*^-}  |R_{\lambda_j}v(t)|^2\le C_6 E_2^2(\bar t).
\eeqs

Then letting $\bar t\to T_*^-$ and using \eqref{limE2} give
\beq\label{Rlamv1}
\lim_{t\to T_*^-}  |R_{\lambda_j}v(t)|^2=0.
\eeq

\medskip
We estimate $|(I_n-R_{\Lambda})v(t)|$ now.
We have 
\beq\label{InR} 
|(I_n-R_{\Lambda})v(t)|
\le \sum_{1\le j\le d, \lambda_j>  \Lambda } |R_{\lambda_j}v(t)|
+\sum_{1\le j\le d, \lambda_j<  \Lambda } |R_{\lambda_j}v(t)|.
\eeq

Taking $t\to T_*^-$ and using \eqref{Rv2} for the first sum and  \eqref{Rlamv1} for the second sum, we obtain the limit  \eqref{remvstar}.

In the case $\Lambda=\Lambda_1$, the second summation in \eqref{InR} is void, hence, \eqref{rvstar1} follows \eqref{Rv2}.

Consider the case $\Lambda>\Lambda_1$. 
To prove \eqref{rvstar2}, we derive more explicit rate of convergence in \eqref{Rlamv1}. We  estimate the right hand side of \eqref{RE3} further.
With $T_*-t$ sufficiently small, hence, less than $1$, and $\varep\in(0,1)$, taking $\bar t=\bar t_\varep\eqdef T_*-(T_*-t)^\varep$ in \eqref{RE3} yields
\beq\label{RE4}
  |R_{\lambda_j}v(t)|^2
\le (T_*-t)^{(1-\varep)C_3\mu}+C_6 E_2^2(\bar t_\varep).
\eeq

Then using \eqref{Rv2} for the first sum  in \eqref{InR}, and  \eqref{RE4} for the second sum  in \eqref{InR}, we obtain
\beq\label{rvstar3}
|(I_n-R_{\Lambda})v(t)|=\bigo((T_*-t)^{(1-\varep)\theta_0}+E_2(\bar t_\varep)+E_2(t)),
 \eeq
where  
\beq\label{thz}
 \theta_0=C_3\mu/2.
 \eeq
 
 Note that $\bar t_\varep<t$ which implies
 \beq \label{Etep}
 E_2(\bar t_\varep)\ge E_2(t).
 \eeq 
 Therefore, we obtain \eqref{rvstar2}  from \eqref{rvstar3}.
 \end{proof}

\begin{remark}\label{pfrm}
The following remarks are in order.
\begin{enumerate}[label=\rnum]
\item By taking $t,\bar t$ sufficiently close to $T_*$ and using \eqref{ytclose} intead of \eqref{yplus} in estimating $H(y(\tau))$ in \eqref{simR}, we can replace $C_3$ in \eqref{RE3} with $c_1( \alpha(c_2\Lambda_n+\varep_0) )^{-1}$. Therefore, we can explicitly replace $\theta_0$ given in \eqref{thz} with any number $\theta_0$ in $(0, c_1\mu /(2\alpha c_2\Lambda_n))$. Note that, unlike \eqref{thz},  this new choice of $\theta_0$ is independent of the solution $y(t)$.

\item\label{rmkb} When $\Lambda_1<\Lambda_n$, 
either  estimate \eqref{rvstar1} or \eqref{rvstar2} for $V_1(t)$  holds true.
Thanks to \eqref{Etep},  estimate \eqref{rvstar1} implies \eqref{rvstar2} for any $\theta_0>0$.
Therefore, in both cases,  there exists $\theta_0>0$ such that, as $t\to T_*^-$,
\beq \label{V1b}
V_1(t)=\bigo((T_*-t)^{(1-\varep)\theta_0}+E_2(T_*-(T_*-t)^\varep))\text{ for all }\varep\in(0,1).
\eeq 

\item The estimates for $V_1(t)$ by \eqref{rvstar1} and \eqref{rvstar2} may not be sharp.
They can be derived directly and differently depending on a particular case of $E_0(t)$, see estimate (6.26) in \cite{H8} for an example.
\end{enumerate}
\end{remark}

As a consequence of Proposition \ref{lem2},  one has
\beq\label{ebarlim}
\lim_{t\to T_*^-} V_1(t)=0.
\eeq

We will use the quantity $V_1(t)$ to describe the asymptotic behavior of $R_\Lambda v(t)$, as $t\to T_*^-$. 

\begin{assumption}\label{moreE}
In the remainder of this section, we assume
\beq\label{E1ss}
 Z_3\eqdef \int_{t_0}^{T_*}\frac{E_1(\tau)}{T_*-\tau}\d \tau<\infty.
 \eeq 
\end{assumption}

\begin{proposition}\label{lem4}
There exists a unit vector $v_*\in\R^n$ such that
\beq\label{RLv}
|R_{\Lambda}v(t)-v_*|=\bigo(V_1(t)+E_1(t))\text{ as $t\to T_*^-$.}
\eeq
\end{proposition}
\begin{proof}
Starting with the norm $|R_\Lambda v(t)|$, we have
\beq\label{oneR}
\big| 1-|R_{\Lambda}v(t)| \big|
=\big| |v(t)|-|R_{\Lambda}v(t)|\big|\le |v(t)-R_{\Lambda}v(t)|=V_1(t).
\eeq
Thus, 
\beq \label{nR}
\lim_{t\to T_*^-}  |R_\Lambda v(t)|= 1.
\eeq 

Then there is $T_0\in(t_0,T_*)$ such that 
\beq \label{Rhalf}
|R_\Lambda v(t)|\ge 1/2, \text{ and, consequently, } R_{\Lambda}v(t)\ne 0,
\text{ for all $t\in[T_0,T_*)$.}
\eeq 

Applying $R_\Lambda$ to equation \eqref{dv} yields, for  $t\in(t_0,T_*)$,
\beq \label{dRv}
\ddt R_{\Lambda}v
 =H(y)(\Lambda-\lambda) R_{\Lambda}v +R_\Lambda g(t).
\eeq

Then,  for $t\in[T_0,T_*)$,
\beq\label{RLveq}
\ddt |R_{\Lambda}v|
=\frac{1}{|R_{\Lambda}v|} \left(\ddt R_{\Lambda}v\right)\cdot R_{\Lambda}v
 =H(y)(\Lambda-\lambda)| R_{\Lambda}v| +g_1(t),
\eeq
where
\beqs
g_1(t)=\frac{R_\Lambda  g(t)\cdot R_{\Lambda}v(t)}{|R_{\Lambda}v(t)|}.
\eeqs

Solving for solution $ |R_{\Lambda}v(t)|$ by the variation of constants formula from the  differential equation \eqref{RLveq} gives, for $\bar t,t\in[T,T_*)$ with $t>\bar t$, 
\beqs
 |R_{\Lambda}v(t)|=e^{\int_{\bar t}^t H(y(\tau))(\Lambda-\lambda(\tau))\d\tau } \left(| R_{\Lambda}v(\bar t) |
 + \int_{\bar t}^t e^{-\int_{\bar t}^\tau H(y(s))(\Lambda-\lambda(s))\d s}  g_1(\tau) \d\tau \right).
\eeqs

It yields
\beq \label{hl0}
\begin{aligned}
&\int_{\bar t}^t H(y(\tau))(\Lambda-\lambda(\tau))\d\tau\\
&=-\ln \left(| R_{\Lambda}v(\bar t) |
 + \int_{\bar t}^t e^{-\int_{\bar t}^\tau H(y(s))(\Lambda-\lambda(s))\d s}  g_1(\tau) \d\tau \right)
   + \ln |R_{\Lambda}v(t)|.
\end{aligned}
\eeq

We have from \eqref{Rcontract}, \eqref{gMy} and \eqref{yplus} that, for $t\in [T,T_*)$,
\beq\label{gg}
|g_1(t)|\le|R_\Lambda  g(t)| \le  |g(t)|\le 2C_2^{\alpha} (T_*-t)^{-1}E_0(t).
\eeq

By \eqref{llH2}, we have, for $\tau\in[\bar t,T_*)$,
$$\Lambda-\lambda(s)\ge -C_5 E_1(s).$$

Combining this inequality with \eqref{Hyt} and \eqref{E1ss} gives
\beqs 
-\int_{\bar t}^\tau H(y(s))(\Lambda-\lambda(s))\d s 
\le C_4C_5 \int_{\bar t}^{\tau}\frac{E_1(s)}{T_*-s}\d s
\le C_4 C_5 Z_3.
\eeqs

This fact and \eqref{gg} imply
\beqs
e^{-\int_{\bar t}^\tau H(y(s))(\Lambda-\lambda(s))\d s} | g_1(\tau) |
\le   C_7 (T_*-\tau)^{-1}E_0(\tau),
\eeqs
where $C_7=2C_2^\alpha e^{C_4C_5 Z_3}$.
Thanks to this and \eqref{condE}, 
\beq\label{eta0}
\lim_{t\to T_*^-}\int_{\bar t}^t e^{-\int_{\bar t}^\tau H(y(s))(\Lambda-\lambda(s))\d s}  g_1(\tau) \d\tau
=\int_{\bar t}^{T_*} e^{-\int_{\bar t}^\tau H(y(s))(\Lambda-\lambda(s))\d s}  g_1(\tau) \d\tau
=\eta(\bar t)\in \R.
\eeq

Note that 
\beq\label{eta}
|\eta(\bar t)|\le  C_7 \int_{\bar t}^{T_*} (T_*-\tau)^{-1}E_0(\tau)\d\tau
=C_7E_1(\bar t).
\eeq

Passing to the limit as $t\to T_*^-$ in \eqref{hl0} and taking into account \eqref{eta0} and \eqref{nR}, we have
\beq\label{hl1}
\int_{\bar t}^{T_*} H(y(\tau))(\Lambda-\lambda(\tau))\d\tau= -\ln (| R_{\Lambda}v(\bar t) |
 + \eta(\bar t))\in\R.
\eeq

By \eqref{hl1}, we can define, for $t\in[T,T_*)$,
$$h(t)=\int_{t}^{T_*} H(y(\tau))(\Lambda-\lambda(\tau))\d\tau\in\R.$$

We rewrite \eqref{hl1} for $\bar t=t$ as
\beqs
h(t)=-\ln (| R_{\Lambda}v(t) | + \eta(t))=-\ln (1+(| R_{\Lambda}v(t) | -1)+ \eta(t)).
\eeqs

With this expression and properties \eqref{oneR} and \eqref{eta}, we have, as $t\to T_*^-$, 
\beq\label{hl2}
| h(t) |
=\bigo(\big|| R_{\Lambda}v(t) | -1\big|+ |\eta(t)|)
=\bigo(V_1(t)+E_1(t)).
\eeq

Solving for $R_{\Lambda}v(t)$ from \eqref{dRv} by the variation of constants formula, one has
\beq\label{RvH}
R_{\Lambda}v(t)
 =e^{\int_{\bar t}^t H(y(\tau))(\Lambda-\lambda(\tau))\d\tau } \left( R_{\Lambda}v(\bar t) 
 + \int_{\bar t}^t e^{-\int_{\bar t}^\tau H(y(s))(\Lambda-\lambda(s))\d s}  R_{\Lambda}g(\tau)\d\tau \right) .
\eeq

Using the same arguments as those from \eqref{gg} to \eqref{eta} with $R_\Lambda g(\tau)$ replacing $g_1(\tau)$, we obtain, similar to \eqref{eta0} and \eqref{eta}, that
\beqs
\lim_{t\to T_*^-} \int_{\bar t}^{t} e^{-\int_{\bar t}^\tau H(y(s))(\Lambda-\lambda(s))\d s}  R_{\Lambda}g(\tau)\d \tau 
= \int_{\bar t}^{T_*} e^{-\int_{\bar t}^\tau H(y(s))(\Lambda-\lambda(s))\d s}  R_{\Lambda}g(\tau) \d\tau=X(\bar t)\in\R^n
\eeqs
for all $\bar t\in[T,T_*)$, and  
\beq \label{oX}
 |X(\bar t)|=\bigo(V_1(t)+E_1(t)) \text{ as $\bar t\to T_*^-$.}
\eeq

Taking $t\to T_*^-$ in \eqref{RvH} gives
\beq\label{vstar}
\lim_{t\to\infty}R_{\Lambda}v(t)=v_*\eqdef e^{h(\bar t)}( R_{\Lambda}v(\bar t) +X(\bar t))\in\R^n.
\eeq

Because of \eqref{nR} and limit \eqref{vstar}, we have  $|v_*|=1$.
Note that
\beqs
X(t)=\int_{t}^{T_*} e^{h(\tau)-h(t)}  R_{\Lambda}g(\tau) \d\tau.
\eeqs

Using $h(t)$, $X(t)$ and $v_*$, we rewrite \eqref{RvH} as
\beqs
R_{\Lambda}v(t)
=e^{h(\bar t)-h(t)} \left (R_{\Lambda}v(\bar t)+X(\bar t)-\int_t^{T_*} e^{h(\tau)-h(\bar t)} R_{\Lambda}g(\tau) \d\tau\right)
=e^{-h(t)}v_*-X(t).
\eeqs

Thus,
\beqs
|R_{\Lambda}v(t)-v_*|
\le |e^{-h(t)}-1|\cdot |v_*| + |X(t)|.
\eeqs

Using \eqref{hl2} and \eqref{oX},  we deduce, as $t\to T_*^-$,
\beqs
|R_{\Lambda}v(t)-v_*|
=\bigo(|h(t))|+|X(t)|)=\bigo(V_1(t)+E_1(t)).
\eeqs

Therefore, we obtain the desired estimate \eqref{RLv}.
The proof is complete.
\end{proof}

\begin{remark}\label{Econ}
By Fubini's theorem, condition \eqref{E1ss} can be rewritten in terms of $E_0(t)$ as
\beq\label{newE1cond}
 \int_{t_0}^{T_*} \frac{E_0(s)}{T_*-s}\ln\left(\frac{T_*-t_0}{T_*-s}\right)\d s<\infty.
\eeq
In fact, if condition \eqref{newE1cond} is met, then so are both conditions \eqref{condE} and \eqref{E1ss}.
\end{remark}

As a consequence of Propositions \ref{lem2} and \ref{lem4},
\beq\label{Rvlim}
\lim_{t\to T_*^-} R_{\Lambda}v(t)=\lim_{t\to T_*^-} v(t)=v_*,
\eeq 
thus,
\beq \label{evzero}
V_2(t)\eqdef |v(t)-v_*|\to 0\text{ as $t\to T_*^-$.}
\eeq 

More specifically, by \eqref{Ebarate}, \eqref{RLv} and the triangle inequality, one has
\beq\label{Ebound}
V_2(t)=\bigo(V_1(t)+E_1(t)) \text{ as $t\to T_*^-$.}
\eeq

The asymptotic behavior of $y(t)$, as $t\to T_*^-$, requires more information about the function $H$.

\begin{definition}\label{Hsym}
Let $\omega$ be a  nonnegative, increasing function defined on $[0,s_0]$, for some number $s_0>0$, such that  
\beq\label{omlim}  \lim_{s\to 0}\omega(s)=\omega(0)=0
\eeq 
and 
\beq\label{Hosc}
|H(x)-H(v_*)|\le \omega(|x-v_*|)\text{ for all $x\in\mathbb S^{n-1}$ with $|x-v_*|\le s_0$.}
\eeq

With a number $T_1\in [T_0,T_*)$ sufficiently close to $T_*$, define
\beq\label{EEzero}
\mathcal E_0(t)=\omega(V_2(t))+V_1(t)+E_0(t)\text{ for $t\in[T_1,T_*)$.}
\eeq
\end{definition}

It is clear that $T_1$ is chosen such that $V_2(t)\le s_0$ for all $t\in [T_1,T_*)$.

Since $\omega$ is monotone, it is Borel measurable. Because $V_2$ is continuous, we have $\omega\circ V_2$ is (Lebesgue) measurable. One can see that $\omega(V_2(t))$ and $V_1(t)$ are bounded on $[T_1,T_*)$. Together with the fact $E_0\in L^1(t_0,T_*)$, this implies
\beq\label{intab}
\mathcal E_0\in L^1(T_1,T_*).
\eeq

There always exists a function $\omega$ as in Definition \ref{Hsym}. For example, $s_0=1/2$ and 
\beq\label{omega1}
\omega(s)=\max\{|H(x)-H(v_*)|:x\in\mathbb S^{n-1},|x-v_*|\le s\}.
\eeq

\begin{theorem}\label{symthm}
Assume 
\beq \label{tEfin}
\mathcal Z_1\eqdef \int_{T_1}^{T_*} \frac{\mathcal E_0(\tau)}{T_*-\tau}\d \tau < \infty.
\eeq

\begin{enumerate}[label=\tnum]
\item Then
\beq\label{ylimxi}
\lim_{t\to T_*^-} (T_*- t)^{1/\alpha}y(t)=\xi_*\in\R^n\setminus\{0\}.
\eeq
Moreover, $\xi_*$ is an eigenvector of $A$ corresponding to $\Lambda$ satisfying 
\beq\label{xiHA}
\alpha \Lambda H(\xi_*)=1.
\eeq

\item The rate of convergence in \eqref{ylimxi} can be specified as follows.
Define
\beq\label{EEone}
\mathcal E_1(t)=\int_t^{T_*}\frac{\mathcal E_0(\tau)}{T_*-\tau} \d\tau\text{ for $t\in[T_1,T_*)$.}
\eeq
Then one has,  as $t\to T_*^-$,
\beq\label{remy}
(T_*-t)^{1/\alpha}|(I_n-R_\Lambda)y(t)|=\bigo(V_1(t)),
\eeq
\beq\label{Ryxi}
|(T_*-t)^{1/\alpha}R_\Lambda y(t)-\xi_* |=\bigo(\mathcal E_1(t))
\eeq
  and, consequently, 
\beq\label{yxio}
|(T_*- t)^{1/\alpha}y(t)-\xi_*| = \bigo(V_1(t)+\mathcal E_1(t)).
\eeq
\end{enumerate}
\end{theorem}
\begin{proof}
Write $(I_n-R_\Lambda)y(t)=|y(t)|({\rm Id}-R_\Lambda)v(t)$ and using \eqref{yplus} and \eqref{Ebarate}, we obtain \eqref{remy}.

We prove \eqref{Ryxi} next. Note from \eqref{tEfin} that
\beq\label{EElim}
\lim_{t\to T_*^-} \mathcal E_1(t)=0.
\eeq

Applying $R_\Lambda$ to equation \eqref{mainode}, we have
\beq\label{Ry1}
\ddt R_\Lambda y =\Lambda  H(y)R_\Lambda y+R_\Lambda f(t).
\eeq

We rewrite $H(y)$ on the right-hand side of \eqref{Ry1} as
\beqs
H(y(t))=|y(t)|^{\alpha} H(v(t))=|R_\Lambda y(t)|^{\alpha} H(v_*)+g_0(t),
\eeqs
where
\begin{align*}
g_0(t)&=|y(t)|^{\alpha} (H(v(t))-H(v_*))+(|y(t)|^{\alpha} -|R_\Lambda y(t)|^{\alpha})H(v_*)\\
&=|y(t)|^{\alpha} \big\{ H(v(t))-H(v_*)+(1 -|R_\Lambda v(t)|^{\alpha})H(v_*)\big\}.
\end{align*}

Then 
\beq\label{Rygood}
\ddt R_\Lambda y 
= \Lambda H(v_*) |R_\Lambda y|^{\alpha}R_\Lambda y+f_0(t),
\eeq
where
$f_0(t)=\Lambda g_0(t)R_\Lambda y(t) +R_\Lambda f(t)$.

We estimate $|g_0(t)|$ first. As $t\to T_*^-$, by \eqref{Hosc} with $x=v(t)$ and \eqref{evzero}, one has
\beq\label{Hvv}
|H(v(t))-H(v_*)|=\bigo(\omega(|v(t)-v_*|)) =\bigo(\omega(V_2(t))).
\eeq

Since $s\in (1/4,3/4)\mapsto s^\alpha$ is a $C^1$-function, by taking $s=|R_\Lambda v(t)|$,  which goes to $1$ as $t\to T_*^-$, and using estimate  \eqref{oneR}, we derive
\beq\label{1Ral}
\big|1 -|R_\Lambda v(t)|^{\alpha}\big|
=\bigo\left(\big |1-|R_\Lambda v(t)|\big| \right)=\bigo(V_1(t))  \text{ as $t\to T_*^-$.}
\eeq

Combining \eqref{Hvv}, \eqref{1Ral} with \eqref{yplus}, we obtain, as $t\to T_*^-$,
\beq\label{g1}
|g_0(t)|=\bigo\left(|y(t)|^{\alpha} [\omega(V_2(t))+V_1(t)]\right)=\bigo\left((T_*-t)^{-1}[\omega(V_2(t))+V_1(t)]\right).
\eeq

We estimate $|f_0(t)|$ next. As $t\to T_*^-$, we have from  \eqref{yplus} and \eqref{fE}  that
\beq\label{Ry5}
|R_\Lambda y(t)|=\bigo((T_*-t)^{-1/\alpha})\text{ and }
|R_\Lambda f(t)|=\bigo((T_*-t)^{-1/\alpha-1}E_0(t)).
\eeq

Combining \eqref{g1} and \eqref{Ry5} gives, as $t\to T_*^-$, 
\beqs
|f_0(t)|=\bigo( (T_*-t)^{-1/\alpha-1}(\omega(V_2(t))+V_1(t)+E_0(t)) )=\bigo((T_*-t)^{-1/\alpha-1}\mathcal E_0(t)).
\eeqs

Note  from \eqref{Rhalf} and the lower bound of $|y(t)|$ in \eqref{yplus} that
\beqs
|R_\Lambda y(t)|=|y(t)||R_\Lambda v(t)|\ge \frac12|y(t)|\ge  \frac{C_1}2  (T_*-t)^{-1/\alpha}
\text{ for all $t\in[T_0,T_*)$.}
\eeqs 

Therefore, there are  $t_0'\in[T_0,T_*)$  and $M_0>0$ such that
\beqs
|f_0(t)|\le M_0|R_\Lambda y(t)|^{1+\alpha} \mathcal E_0(t) \text{ for all  $t\in[t_0',T_*)$.  } 
\eeqs

We apply Theorem \ref{simthm} to solution $R_\Lambda y(t)$ of equation \eqref{Rygood} on the interval $[t_0',T_*)$. Specifically, $y(t):=R_\Lambda y(t)$ satisfies equation \eqref{basicf} on $(t_0',T_*)$ with constant $a:=\Lambda H(v_*) $ and $f:=f_0$, $E_0:=M_0\mathcal E_0$. 
By \eqref{Rhalf} and \eqref{yplus}, one has $|R_\Lambda y(t)|\ge |y(t)|/2 \to\infty$ as $t\to T_*^-$.
Then there exists a nonzero vector $\xi_*\in\R^n$ such that 
\beq\label{Rxi3}
|(T_*-t)^{1/\alpha}R_\Lambda y(t)-\xi_* |=\bigo(M_0\mathcal E_1(t))
\eeq
and 
\beq\label{Hv1}
|\xi_*|=(\alpha \Lambda H(v_*) )^{-1/\alpha} .
\eeq
The  statement  \eqref{Ryxi} follows \eqref{Rxi3}. Then \eqref{yxio} follows \eqref{remy}, \eqref{Ryxi} and the triangle inequality. Clearly, \eqref{ylimxi} follows \eqref{yxio} and properties \eqref{ebarlim}, \eqref{EElim}.

Because of the limit \eqref{ylimxi} and the fact $\xi_*\ne 0$, we have $\xi_*\in R_\Lambda(\R^n)\setminus\{0\}$. Hence, $\xi_*$ is an eigenvector of $A$ associated with $\Lambda$.

Let $w(t)=(T_*- t)^{1/\alpha}y(t)$ and write 
$v(t)=w(t)/|w(t)|$.
Passing $t\to T_*^-$ and noticing that $v(t)\to v_*$ and $w(t)\to \xi_*$, thanks to \eqref{Rvlim} and \eqref{ylimxi},
we obtain
\beq\label{vxi} 
v_*=\xi_*/|\xi_*|.
\eeq 

Then it follows \eqref{Hv1}, the fact $H$ is positively homegeneous of degree $\alpha$, and relation \eqref{vxi}   that
$$1=\alpha \Lambda H(v_*)  |\xi_*|^{\alpha}=\alpha \Lambda H(|\xi_*|v_* )
=\alpha \Lambda H(\xi_*).$$
Hence, we obtain \eqref{xiHA} and the proof is complete.
\end{proof}

\section{The non-symmetric matrix case}\label{nonsymcase}
Consider the case $A$ is not symmetric.  Let $y(t)$ be a solution of equation \eqref{mainode}  as in the beginning of Section \ref{ratesec}.  
\begin{itemize}
\item Assume the function $E_0(t)$ in \eqref{fE} satisfies \eqref{newE1cond} which, according to Remark \ref{Econ}, is equivalent to \eqref{condE}  and \eqref{E1ss} with $E_1(t)$ being defined by \eqref{E1def}.

\item When $\Lambda_1<\Lambda_n$, we impose Assumption \ref{E2assum} and define $E_2(t)$ as in \eqref{E2def}.
\end{itemize}

Recall that matrices  $A_0$ and $S$ are from \eqref{Adiag}.
For $t\in[t_0,T_*)$, define $\widetilde E_0(t)$ by  \eqref{tilEzero},
\begin{align*}
\widetilde E_1(t)&=\int_{t}^{T_*} (T_*-\tau)^{-1}\widetilde E_0(\tau)\d\tau
\intertext{and, when $\Lambda_1<\Lambda_n$,}
\widetilde E_2(t)&=\left(\int_{t}^{T_*} (T_*-\tau)^{-1}\widetilde E_0^2(\tau)\d\tau\right)^{1/2}.
\end{align*} 

Thanks to  \eqref{tilEzero}, we immediately have 
\beq\label{propEE}
\widetilde E_1(t)=\|S\|\cdot  \|S^{-1}\|^{\alpha+1}E_1(t) \text{ and }
\widetilde E_2(t)=\|S\|\cdot  \|S^{-1}\|^{\alpha+1}E_2(t). 
\eeq 

As a consequence of \eqref{condE}, Assumption \ref{E2assum} and \eqref{propEE}, one has 
\beqs
\widetilde Z_1\eqdef \widetilde E_1(t_0)=\|S\|\cdot  \|S^{-1}\|^{\alpha+1}Z_1<\infty \text{ which also implies } \widetilde E_0\in L^1(t_0,T_*),
\eeqs 
and, when $\Lambda_1<\Lambda_n$, 
\beqs
\widetilde Z_2\eqdef (\widetilde E_2(t_0))^2=\|S\|^2 \|S^{-1}\|^{2(\alpha+1)} Z_2<\infty. 
\eeqs 

From \eqref{propEE}, \eqref{Elim}, \eqref{limE2} and \eqref{E1ss}, it follows that 
\beqs
\lim_{t\to T_*^-} \widetilde E_1(t)=0,\quad 
\lim_{t\to T_*^-} \widetilde E_2(t)=0
\eeqs
and
\beqs
\widetilde Z_3\eqdef \int_{t_0}^{T_*}\frac{\widetilde E_1(\tau)}{T_*-\tau}\d \tau
=\|S\|\cdot  \|S^{-1}\|^{\alpha+1}Z_3<\infty.
 \eeqs

We follow Part 2 of the proof of Lemma \ref{estthm} up to \eqref{tilEzero}.

For $t\in[t_0,T_*)$, recalling $\widetilde y(t)=Sy(t)$, we set $\widetilde v(t)=\widetilde y(t)/|\widetilde y(t)|$. 
The asymptotic behavior of $\widetilde v(t)$ is described below.

\begin{lemma}\label{nonsymlem}
The following statements hold true.
\begin{enumerate}[label=\tnum]
\item There are unique eigenvalue $\Lambda$ of $A$ and unit vector $\widetilde v_*$ such that
\beq\label{EEbar}
\widetilde V_1(t)\eqdef |(I_n-\widehat R_\Lambda)\widetilde v(t)|\to 0 \text{ as $t\to T_*^-$,}
\eeq
and
\beq\label{REE}
|\widehat R_\Lambda \widetilde v(t)-\widetilde v_*|=\bigo(\widetilde V_1(t)+\widetilde E_1(t)).
\eeq
Consequently,
\beq\label{tilV2}
\widetilde V_2(t)\eqdef |\widetilde v(t)-\widetilde v_*|=\bigo(\widetilde V_1(t)+\widetilde E_1(t)) \text{ as $t\to T_*^-$.}
\eeq

\item The $\widetilde V_1(t)$ can be estimated more explicitly as follows.

If $\Lambda_1=\Lambda_n$, then 
\beq \label{EE0}
\widetilde V_1(t)=0\text{  for all $t\in[t_0,T_*)$.}
\eeq

 If $\Lambda_1=\Lambda<\Lambda_n$, then
\beq\label{EE2}
\widetilde V_1(t)=\bigo(\widetilde E_2(t)) \text{ as $t\to T_*^-$.}
 \eeq
 
 If  $\Lambda_1<\Lambda\le\Lambda_n$, then there is $\widetilde\theta_0>0$ such that,  as $t\to T_*^-$,
\beq\label{EeE}
\widetilde V_1(t)=\bigo((T_*-t)^{(1-\varep)\widetilde\theta_0}+\widetilde E_2(T_*-(T_*-t)^\varep))
 \text{ for any $\varep\in(0,1)$.}
 \eeq
\end{enumerate}
\end{lemma}
\begin{proof}
Part (i). We apply Propositions \ref{lem2} and \ref{lem4} to the solution $\widetilde y(t)$ of equation \eqref{zeq}, noting that 
$A_0$ replaces $A$, 
$\widehat R_{\lambda_j}$  replaces $R_{\lambda_j}$, 
$\widetilde H$ and $\widetilde f$ in \eqref{HRz} replace $H$ and $f$ respectively,
$\widetilde E_j(t)$ replaces $E_j(t)$ for $j=0,1,2$,
and 
$\widetilde V_j(t)$ replaces $V_j(t)$ for $j=1,2$.
Then there exist a unique eigenvalue $\Lambda$ of $A_0$ and a unique unit vector $\widetilde v_*$ such that one has \eqref{EEbar} and \eqref{REE}, by \eqref{remvstar} and \eqref{RLv}, respectively.
Then estimate \eqref{tilV2} follows \eqref{EEbar}, \eqref{REE} and the triangle inequality.
Clearly, $\Lambda$ is an eigenvalue $\Lambda$ of $A$.

Part (ii). We have estimates \eqref{EE0}, \eqref{EE2} and  \eqref{EeE} follow \eqref{rvstar0}, \eqref{rvstar1} and \eqref{rvstar2}, respectively.
\end{proof}

Set $v_*=S^{-1} \widetilde v_*$.  Similar to Definition \ref{Hsym}, let $\omega$ be a  nonnegative, increasing function defined on $[0,s_0]$, for some number $s_0>0$ that satisfies \eqref{omlim} and  
\beq\label{Hosc2}
|H(x)-H(v_*)|\le \omega(|x-v_*|)\text{ for all $x\in\mathbb \R^n\setminus\{0\}$ with $|x-v_*|\le s_0$.}
\eeq

Again, such a function $\omega$ exists, for e.g., $s_0=|v_*|/2$ and $\omega(s)$ is as in \eqref{omega1} with $x\in \R^n\setminus\{0\}$. 

Select a number $T_1\in [t_0,T_*)$ sufficiently close to $T_*$ such that 
$\|S^{-1}\|  \widetilde V_2(t)\le s_0$ for all $t\in[T_1,T_*)$, and define
\beq\label{EEorig}
\widetilde{\mathcal E}_0(t)=\omega(\|S^{-1}\| \widetilde V_2(t))+\widetilde V_1(t)+\widetilde E_0(t)
\text{ for $t\in[T_1,T_*)$.}
\eeq

Note that the function $\omega(\|S^{-1}\| \widetilde V_2(t))$ in \eqref{EEorig} is measurable.
Similar to \eqref{intab}, one has
$\widetilde{\mathcal E}_0\in L^1(T_1,T_*)$.

We are ready to establish the asymptotic behavior of $y(t)$ for $t$ near $T_*$.

\begin{theorem}\label{nonsymthm}
Assume 
\beq\label{cEtil1}
\widetilde {\mathcal Z}_1\eqdef 
\int_{T_1}^{T_*}\frac{\widetilde{\mathcal E}_0(\tau)}{T_*-\tau} \d\tau<\infty
\eeq
and define
\beq\label{defEtil1}
\widetilde{\mathcal E}_1(t)=\int_t^{T_*}\frac{\widetilde{\mathcal E}_0(\tau)}{T_*-\tau} \d\tau \text{ for $t\in[T_1,T_*)$.}
\eeq

Then the statements in part {\rm (i)} of Theorem \ref{symthm} hold true.
More specifically, one has, as $t\to T_*^-$,  
\beq\label{symo1}
(T_*-t)^{1/\alpha}|(I_n-R_\Lambda) y(t)|=\bigo(\widetilde V_1(t)),
\eeq
\beq\label{symo2}
|(T_*-t)^{1/\alpha}R_\Lambda y(t)-\xi_*|=\bigo(\widetilde{\mathcal E}_1(t))
\eeq
and, consequently,
\beq\label{symo3}
|(T_*-t)^{1/\alpha}y(t)-\xi_*|=\bigo(\widetilde V_1(t)+\widetilde{\mathcal E}_1(t)).
\eeq
\end{theorem}
\begin{proof}
First, we notice from \eqref{cEtil1} and  \eqref{defEtil1} that
\beq\label{EE1lim}
\lim_{t\to T_*^-} \widetilde{\mathcal E}_1(t)=0.
\eeq

Define 
\beq \label{oomg}
\widetilde \omega(s)=\omega(\|S^{-1}\|s)\text{ for $s\in[0,\widetilde s_0]$, where $\widetilde s_0=s_0/\|S^{-1}\|$.}
\eeq 

Then $\widetilde\omega$ is a nonnegative, increasing function on $[0,\widetilde s_0]$ and  
$\lim_{s\to 0} \widetilde\omega(s)=\widetilde\omega(0)=0$. Again, the monotonicity of $\widetilde\omega$ implies its Lebesgue-measurability.

Suppose $z\in\R^n$ with $|z-\widetilde v_*|\le \widetilde s_0$, then  we have $|S^{-1}z-v_*|=|S^{-1}(z-\widetilde v_*)|\le s_0$. Thus,
\begin{align*}
|\widetilde H(z)-\widetilde H(\widetilde v_*)|
&=|H(S^{-1}z)-H(v_*)|  \le \omega(|S^{-1}z-v_*|)\\
&\le \omega(\|S^{-1}\|\cdot |z-\widetilde v_*|)=\widetilde\omega(|z-\widetilde v_*|).
\end{align*}

The $\widetilde{\mathcal E}_0(t)$ in \eqref{EEorig}, in fact, is
\beq\label{EEoo}
\widetilde{\mathcal E}_0(t)=\widetilde\omega(\widetilde V_2(t))+\widetilde V_1(t)+\widetilde E_0(t) 
\eeq 
which resembles \eqref{EEzero}.

By Theorem \ref{symthm} applied to solution $\widetilde y(t)$ of equation \eqref{zeq}, there exists an eigenvector $\widetilde \xi_*$ of $A_0$ associated  with $\Lambda$ such that 
\beq\label{hatRz}
(T_*-t)^{1/\alpha}|(I_n-\widehat R_\Lambda)\widetilde y(t)|=\bigo(\widetilde V_1(t)),\quad
|(T_*-t)^{1/\alpha}\widehat R_\Lambda \widetilde y(t)-\widetilde \xi_* |=\bigo(\widetilde{\mathcal E}_1(t))
\eeq
as $t\to T_*^-$, and
\beq\label{A0xi} 
\alpha \Lambda \widetilde H(\widetilde \xi_*)=1.
\eeq

Let $\xi_*=S^{-1}\widetilde \xi_*$. It is elementary to see that $\xi_*$ is an eigenvector of $A$ associated  with $\Lambda$. Also,   \eqref{xiHA} follows  \eqref{A0xi} and the relation $\widetilde H(\widetilde \xi_*)=H(\xi_*)$.

Rewriting \eqref{hatRz} in terms of  $R_\Lambda$ and $y(t)$, we obtain 
\beqs
(T_*-t)^{1/\alpha}|S(I_n-R_\Lambda)y(t)|=\bigo(\widetilde V_1(t)),\quad
|S( (T_*-t)^{1/\alpha}R_\Lambda y(t)-\xi_* )|=\bigo(\widetilde{\mathcal E}_1(t)).
\eeqs

These estimates imply \eqref{symo1} and \eqref{symo2}, and, with the use of the triangle inequality, \eqref{symo3} follows.
Finally, \eqref{EEbar}, \eqref{EE1lim} and \eqref{symo3} imply \eqref{ylimxi}.   
\end{proof}

\section{Some typical results}\label{egsec}

Regarding the function $H$, we specify the following class of functions.

\begin{definition}\label{HCdef}
Let $E$ be a nonempty subset of $\R^n$ and $F$ be a function from $E$ to $\R$. We say $F$ has property (HC) on $E$ if, 
for any  $x_0\in E$,  there exist  numbers  $r,C,\gamma>0$ such that
\beq\label{FHder}
|F(x)-F(x_0)|\le C|x-x_0|^\gamma
\eeq
for any $x\in E$ with $|x-x_0|<r$.
\end{definition}

 For example, if $E$ is open and $F$ is locally H\"older continuous of order $\beta\in(0,1]$ on $E$, then $F$ has property (HC) on $E$ with $\gamma=\beta$ in \eqref{FHder} for all $x_0\in E$.
 For more specific examples for the function $H$, see  \cite[Example 3.4]{H8}, \cite[Example 5.7]{H7} and \cite[Section 6]{CaHK1}.

\begin{lemma}\label{Hcts}
Let $F\in \mathcal H_{\alpha}(\R^n,\R)$ for some $\alpha>0$.
\begin{enumerate}[label=\tnum]
\item\label{HH1} If $F>0$ on $\mathbb S^{n-1}$, then $F>0$ on  $\R^n\setminus \{0\}$.

\item\label{HH2} If $F$ is continuous on  $\mathbb S^{n-1}$, then  it is continuous on $\R^n\setminus \{0\}$.

Assume $F$ has property (HC) on  $\mathbb S^{n-1}$ in \ref{HH3}--\ref{HH5} below.

\item\label{HH3} Then $F$ has property (HC) on $\R^n\setminus\{0\}$.

\item\label{HH4} If  $\varphi$ is a function from $\R^n\setminus\{0\}$ to $\R^n\setminus\{0\}$ that has property (HC) on $\R^n\setminus\{0\}$. Then $F\circ \varphi$ has property (HC) on $\R^n\setminus \{0\}$.

\item\label{HH5}  If $K$ is an invertible $n\times n$ matrix, then the function $x\in \R^n\setminus\{0\} \mapsto F(Kx)$ has property (HC) on $\R^n\setminus \{0\}$.
\end{enumerate}
\end{lemma}
\begin{proof}
This lemma is proved in \cite[Lemma 4.1]{H8} for positively homogeneous functions of negative degrees, but it holds true also for positive degrees as stated above.

To track the power $\gamma$ in \eqref{FHder} explicitly, we present full calculations for part \ref{HH5} here.
Let $x_0\in \R^n\setminus\{0\}$. Take any $x\in\R^n$ sufficiently close to $x_0$. Then $x\ne 0$, and $Kx/|Kx|\in \mathbb S^{n-1}$ is close to $Kx_0/|Kx_0|\in \mathbb S^{n-1}$. Because $F$ has property (HC) on  $\mathbb S^{n-1}$, there exist $C>0$ and $\gamma>0$ such that
\beqs 
\left | F(Kx/|Kx|)-F(Kx_0/|Kx_0|)\right| \le C\left| \frac{Kx}{|Kx|}-\frac{Kx_0}{|Kx_0|}\right|^\gamma.
\eeqs 

Below, $C'$ denotes a generic positive constant.
We have
\begin{align*}
 &|F(Kx)-F(Kx_0)|
 =  \Big| |Kx|^{\alpha}F(Kx/|Kx|)-|Kx_0|^{\alpha}F(Kx_0/|Kx_0|)\Big|\\
 &\le |Kx|^{\alpha}\Big | F(Kx/|Kx|)-F(Kx_0/|Kx_0|)\Big| 
 +\Big | |Kx|^{\alpha}- |Kx_0|^{\alpha}\Big| \cdot\big|F(Kx_0/|Kx_0|)\big| \\
 &\le C' \left| \frac{Kx}{|Kx|}-\frac{Kx_0}{|Kx_0|}\right|^\gamma +C'\Big | |Kx|^{\alpha}- |Kx_0|^{\alpha}\Big|.
\end{align*}

Note  that functions $x\in E\mapsto Kx/|Kx|$ and $x\in E\mapsto |Kx|^{\alpha}$ are $C^1$-functions.
Hence, we obtain, for $x$ sufficiently close to $x_0$,
\beq\label{FKHC} 
 |F(Kx)-F(Kx_0)|
\le C' | x-x_0|^\gamma +C'|x-x_0|\\
 \le C'|x-x_0|^{\min\{1,\gamma\}},
\eeq
which proves \ref{HH5}.
\end{proof}

\medskip
\noindent\textbf{Part I.}
We consider a solution $y(t)$ of equation \eqref{mainode} as in the beginning of Section \ref{ratesec}, and  follow the notation in Sections \ref{symcase} and \ref{nonsymcase}. The big oh in calculations below is for $t\to T_*^-$. 

\begin{theorem}\label{eg1}
Assume the function $H$  has property (HC) on  $\mathbb S^{n-1}$, 
and, in \eqref{fE},
\beq \label{epower}
E_0(t)=M(T_*-t)^\delta\text{ for some numbers $M,\delta>0$.}
\eeq 

Then there exists an eigenvector $\xi_*$  of $A$ associated with an eigenvalue $\Lambda$ such that \eqref{xiHA} holds and, as $t\to T_*^-$, 
\beq \label{sc1}
|  (T_*-t)^{1/\alpha}y(t)-\xi_*|=\bigo((T_*-t)^{\varep_*}) \text{ for some number $\varep_*>0$.}
\eeq
\end{theorem}
\begin{proof}
The calculations below will allow us to apply Theorems \ref{symthm} and \ref{nonsymthm}.
Therefore,  there exists an eigenvector $\xi_*$  of $A$ associated with an eigenvalue $\Lambda$ such that \eqref{xiHA}  and the limit \eqref{ylimxi} hold.
We will clarify the estimates \eqref{yxio} and \eqref{symo3}.

\medskip
\noindent\textit{Case $A$ is symmetric.}
By \eqref{E1def},
\beqs
E_1(t)=\frac{M(T_*-t)^\delta}{\delta}.
\eeqs

It follows that $Z_1=E_1(t_0)<\infty$  and 
\beqs
Z_3=\int_{t_0}^{T_*} \frac{E_1(t)}{T_*-t}\d t=\frac{M(T_*-t_0)^\delta}{\delta^2}<\infty.
\eeqs
Therefore, conditions  \eqref{condE} and  \eqref{E1ss} are met.

Since  $H$  has property (HC) on  $\mathbb S^{n-1}$, there are some numbers $s_0,C,\gamma >0$ such that the function $\omega$ in Definition \ref{Hsym} is
\beq \label{opower}
\omega(s)=Cs^\gamma \text{ for all $s\in[0,s_0]$.}
\eeq 

\begin{itemize}
\item Case $\Lambda_1=\Lambda_n$. By \eqref{rvstar0} and \eqref{Ebound},
\beqs 
V_2(t)=\bigo((T_*-t)^\delta).
\eeqs
By  \eqref{rvstar0}, \eqref{EEzero} and \eqref{EEone},
\beq\label{me1}
\mathcal E_0(t)=\bigo((T_*-t)^{\varep_1}) \text{  and }
\mathcal E_1(t)=\bigo((T_*-t)^{\varep_1}),
\eeq
where $\varep_1=\min\{\gamma\delta,\delta\}$.
Therefore, following \eqref{yxio},
\beq\label{conv0}
|(T_*- t)^{1/\alpha}y(t)-\xi_*| = \bigo((T_*-t)^{\varep_1}).
\eeq

\item Case $\Lambda_1<\Lambda_n$. By  \eqref{E2def},
\beqs
E_2(t)=\frac{M(T_*-t)^\delta}{\sqrt{2\delta}}. 
\eeqs
Then $Z_2=(E_2(t_0))^2<\infty$ which yields \eqref{condE2}.
Let $\varep\in(0,1)$. Then
\beqs E_2(T_*-(T_*-t)^\varep)=\frac{M(T_*-t)^{\varep\delta}}{\sqrt{2\delta}}.
\eeqs  
From \eqref{V1b} and \eqref{Ebound}, one has
\beq\label{VV} 
V_1(t)=\bigo((T_*-t)^{\varep_2}), \quad
V_2(t)=\bigo((T_*-t)^{\varep_2}), 
\eeq
where $\varep_2=\min\{(1-\varep)\theta_0,\varep\delta\}\in (0,\delta)$. 
By \eqref{EEzero}, \eqref{opower} and \eqref{EEone},
\beq\label{doubleE}
\mathcal E_0(t)=\bigo((T_*-t)^{\varep_3}) \text{  and }
\mathcal E_1(t)=\bigo((T_*-t)^{\varep_3}),
\eeq
where $\varep_3=\min\{\gamma\varep_2,\varep_2,\delta\}=\varep_2\min\{1,\gamma\}\le \varep_2$.
Therefore, following \eqref{yxio},
\beq\label{conv1}
|(T_*- t)^{1/\alpha}y(t)-\xi_*| = \bigo((T_*-t)^{\varep_3}).
\eeq
\end{itemize} 

\medskip
\noindent\textit{Case $A$ is non-symmetric.} 
By \eqref{tilEzero} and \eqref{epower}, we have
\beqs
\widetilde E_0(t)=\|S\| \cdot \|S^{-1}\|^{\alpha+1}M(T_*-t)^\delta
\eeqs
which has the form of \eqref{epower}.

Thanks to Lemma \ref{Hcts}\ref{HH3}, the function $H$ has property (HC) on $\R^n\setminus\{0\}$. Then the function $\omega$ in \eqref{Hosc2} can be taken of the form \eqref{opower}. Thus, the function $\widetilde \omega$ in \eqref{oomg} is 
\beqs
\widetilde \omega(s)=C\| S^{-1}\|^\gamma s^\gamma,
\eeqs
which has the same form as \eqref{opower}. 
Hence, we can apply the above calculations for the case of symmetric $A$ to the solution $\widetilde y(t)$.

\begin{itemize}
\item Case $\Lambda_1=\Lambda_n$. From \eqref{me1}, we have
\beqs
\widetilde{\mathcal E}_1(t)=\bigo((T_*-t)^{\varep_1}).
\eeqs
Thanks to this, \eqref{EE0} and \eqref{symo3}, we obtain \eqref{conv0} again.
 
\item Case $\Lambda_1<\Lambda_n$. 
 Similar to \eqref{VV} and \eqref{doubleE}, one has
\beqs  \widetilde V_1(t)=\bigo((T_*-t)^{\varep_2})
\text{ and }
\widetilde {\mathcal E}_1(t)=\bigo((T_*-t)^{\varep_3}).
\eeqs
Combining this with \eqref{symo3} and the fact $\varep_3\le\varep_2$, we obtain \eqref{conv1} again.
\end{itemize}

\medskip
\noindent\textit{Summary.} By estimates \eqref{conv0} and \eqref{conv1}  for both cases of symmetric and non-symmetric $A$, we conclude that there is $\varep_*>0$ such that \eqref{sc1} holds true.
\end{proof}

\begin{theorem}\label{eg2}
Assume the function $H$  has property (HC) on  $\mathbb S^{n-1}$ with the same power $\gamma\in(0,1]$ in \eqref{FHder} for $F=H$ and all $x_0\in\mathbb S^{n-1}$. Suppose  $E_0(t)$  in \eqref{fE} is
\beq\label{elog}
E_0(t)=\frac{M}{|\ln (T_*-t)|^p} \text{ for all $t\in[t_0,T_*)$ sufficiently close to $T_*$,}
\eeq 
where  $M>0$ and 
\beq\label{pcond}
p>1+1/\gamma.
\eeq

Then there exists an eigenvector $\xi_*$  of $A$ associated with an eigenvalue $\Lambda$ such that \eqref{xiHA} holds and, as $t\to T_*^-$, 
\beq \label{sc2}
|  (T_*-t)^{1/\alpha}y(t)-\xi_*|=\bigo(|\ln(T_*-t)|^{-p_*})\text{ with $p_*=\gamma(p-1)-1>0$.}
\eeq
\end{theorem}
\begin{proof}
Without loss of generality, we can assume $T_*-t_0<1$ and \eqref{elog} holds for all $t\in[t_0,T_*)$.
We will apply Theorems \ref{symthm} and \ref{nonsymthm} to obtain the desired $\xi_*$ and $\Lambda$.
Our calculations below will justify such applications and elaborate the estimates \eqref{yxio} and \eqref{symo3}.

\medskip
\noindent\textit{Case $A$ is symmetric.}
By \eqref{E1def},
\beqs
E_1(t)=\frac{M}{(p-1)|\ln(T_*-t)|^{p-1}} \text{ and } Z_1=E_1(t_0)<\infty.
\eeqs
Hence, condition \eqref{condE} is satisfied.
Since $\gamma\in(0,1]$, we  have from \eqref{pcond} that $p>2$.
Thus,
\beqs
Z_3=\int_{t_*}^{T_*} \frac{E_1(t)}{T_*-t}\d t=\frac{M}{(p-1)(p-2)|\ln(T_*-t_*)|^{p-2}}<\infty
\eeqs
which verifies \eqref{E1ss}.

Thanks to the condition on $H$, there exist $s_0,C>0$ such that the function $\omega$ in Definition \ref{Hsym} is the same as in \eqref{opower}.

\begin{itemize}
\item Case $\Lambda_1=\Lambda_n$. By \eqref{rvstar0} and \eqref{Ebound},
\beqs 
V_2(t)=\bigo(1/|\ln(T_*-t)|^{p-1}).
\eeqs
By  \eqref{rvstar0}, \eqref{EEzero} and \eqref{EEone},
\beq\label{he1}
\mathcal E_0(t)=\bigo(1/|\ln(T_*-t)|^{p_1})\text{  and }
\mathcal E_1(t)=\bigo(1/|\ln(T_*-t)|^{p_1-1}),
\eeq
where $p_1=\min\{\gamma(p-1),p\}=\gamma(p-1)$.
Thanks to condition \eqref{pcond}, we have $p_1>1$.
Therefore, following \eqref{yxio},
\beq\label{conv2}
|(T_*- t)^{1/\alpha}y(t)-\xi_*| =\bigo(1/|\ln(T_*-t)|^{p_1-1}).
\eeq

\item Case $\Lambda_1<\Lambda_n$. 
By  \eqref{E2def},
\beqs
E_2(t)=\frac{M}{\sqrt{2p-1}|\ln(T_*-t)|^{p-1/2}}. 
\eeqs
Then, again, $Z_2=(E_2(t_0))^2<\infty$ which gives \eqref{condE2}.
Let $\varep\in(0,1)$. We have
\beqs
E_2(T_*-(T_*-t)^\varep)=\frac{M}{\sqrt{2p-1}\,\varep^{p-1/2} |\ln(T_*-t)|^{p-1/2}}.
\eeqs 
Using  \eqref{V1b} and \eqref{Ebound}, we obtain
\beq\label{he3}
 V_1(t)=\bigo(1/|\ln(T_*-t)|^{p-1/2})
\text{ and }
V_2(t)=\bigo(1/|\ln(T_*-t)|^{p-1}).
\eeq
By \eqref{EEzero},
\beqs
\mathcal E_0(t)=\bigo(1/|\ln(T_*-t)|^{p_2}) \text{ with $p_2=\min\{\gamma(p-1),p-1/2,p\}=p_1$.}
\eeqs
Then, by \eqref{EEone},
\beq \label{he2}
\mathcal E_1(t)=\bigo(1/|\ln(T_*-t)|^{p_1-1}).
\eeq 
Therefore, following \eqref{yxio} and the fact $p_1<p$, we obtain
\beq\label{conv3}
|(T_*- t)^{1/\alpha}y(t)-\xi_*| = \bigo(1/|\ln(T_*-t)|^{p_1-1}).
\eeq
\end{itemize}

\medskip
\noindent\textit{Case $A$ is non-symmetric.} By \eqref{tilEzero} and \eqref{elog},
\beqs
\widetilde E_0(t)=\frac{\|S\| \cdot \|S^{-1}\|^{\alpha+1}M}{|\ln (T_*-t)|^p} ,
\eeqs
which has the form of \eqref{elog}. It follows \eqref{FKHC} with $F=H$ and $K=I_n$ that the function $\omega$ in \eqref{Hosc2} can be chosen to be $\omega(s)=Cs^{\gamma}$ again.

\begin{itemize}
\item Case $\Lambda_1=\Lambda_n$. 
  Following the same calculations as in the case $A$ is symmetric, we have from \eqref{he1} that
\beqs
\widetilde{\mathcal E}_1(t)=\bigo(1/|\ln(T_*-t)|^{p_1-1}).
\eeqs
Thanks to this, \eqref{EE0} and \eqref{symo3}, we obtain
\beq\label{conv22}
|(T_*- t)^{1/\alpha}y(t)-\xi_*| =\bigo(1/|\ln(T_*-t)|^{p_1-1}).
\eeq
 
\item Case $\Lambda_1<\Lambda_n$. 
Following the same calculations as in the case $A$ is symmetric, using \eqref{EEoo}, we have from \eqref{he3} and \eqref{he2} that 
\beqs 
\widetilde V_1(t)=\bigo(1/|\ln(T_*-t)|^{p-1/2})
\text{ and }
\widetilde{\mathcal E}_1(t)=\bigo(1/|\ln(T_*-t)|^{p_1-1}).
\eeqs 
By \eqref{symo3}, we obtain
\beq\label{conv33}
|(T_*-t)^{1/\alpha}y(t)-\xi_*|= \bigo(1/|\ln(T_*-t)|^{p_1-1}).
\eeq
\end{itemize}

\medskip
\noindent\textit{Summary.} By  \eqref{conv2}, \eqref{conv3},  \eqref{conv22} and \eqref{conv33}, we obtain \eqref{sc2}.
\end{proof}

We note that it is always possible to assume  $\gamma\le1$ as stated in Theorem \ref{eg2}. This is a convenience, not a restriction.

\medskip
\noindent\textbf{Part II.}
We turn to solutions of equation \eqref{yGeq}  now.
We consider the following particular cases for the functions $H$ and $G$. Let $t_0\ge 0$.

\medskip
\textit{Case 1.} The function $H$  has property (HC) on  $\mathbb S^{n-1}$,  $G(t,x)$ is a continuous function from $[t_0,\infty)\times (\R^n\setminus\{0\})$ to $\R^n$,   and there are $M,\delta>0$ and $r_*>1$ such that  
\beq \label{G1}
|G(t,x)|\le M|x|^{1+\alpha-\delta} \text{ for all  $t\ge t_0$ and $|x|\ge r_*$.}
\eeq

\medskip
\textit{Case 2.} The function $H$  has property (HC) on  $\mathbb S^{n-1}$ with the same power $\gamma\in(0,1]$ in \eqref{FHder} for $F=H$ and all $x_0\in\mathbb S^{n-1}$, $G(t,x)$ is a continuous function from $[t_0,\infty)\times (\R^n\setminus\{0\})$ to $\R^n$, and there are $M>0$, $r_*>e$ and $p>1+1/\gamma$  such that  
\beq \label{G2}
 |G(t,x)|\le \frac{M|x|^{1+\alpha}} {\left ( \ln|x| \right)^p} \text{ for all  $t\ge t_0$ and $|x|\ge r_*$.}
 \eeq 

It is clear that $G(t,x)$, in both cases,  satisfies \eqref{Grate}.
According to Theorem \ref{blowthm}, there exist finite time blow-ups for solutions of \eqref{yGeq} with sufficiently large initial data. 
Their asymptotic behavior near the blow-up time is the following.

\begin{theorem}\label{typthm}
Consider Cases 1 and 2 above.  Suppose $y\in C^1([t_0,T_*),\R^n\setminus\{0\})$, with a number $T_*>t_0$,  satisfies equation \eqref{yGeq} on $(t_0,T_*)$ and has property \eqref{limbig}.

Then there exists an eigenvector $\xi_*$  of $A$ associated with an eigenvalue $\Lambda$ such that \eqref{xiHA} holds and, as $t\to T_*^-$, 
\begin{align}
\label{mainest1}
|  (T_*-t)^{1/\alpha}y(t)-\xi_*|&=\bigo((T_*-t)^{\varep_*}) &\text{ in  Case 1,}\\
\label{mainest2}
|  (T_*-t)^{1/\alpha}y(t)-\xi_*|&=\bigo(|\ln(T_*-t)|^{-p_*}) &\text{ in  Case 2,}
\end{align}
for some numbers $\varep_*,p_*>0$.
\end{theorem}
\begin{proof}
Observe that $y(t)$ is a solution of equation \eqref{mainode}  with $f(t)=G(t,y(t))$ for $t\in[t_0,T_*)$. 

By the virtue of Corollary \ref{boundcor}, we have the  estimates in \eqref{yplus} for $|y(t)|$.
With \eqref{yplus}, we derive more specific estimates for $f(t)$ in order to apply Theorems \ref{eg1}  and \ref{eg2}.

\medskip
In Case 1, by \eqref{yplus} and \eqref{G1}, we have, as  $t<T_*$  close to $T_*$,
\beqs
|f(t)|\le M|y(t)|^{1+\alpha} \cdot C_1^{-\delta}(T_*-t)^{\delta/\alpha}. 
\eeqs
Hence, by considering the solution in a subinterval $[t_0',T_*)$ with $t_0'$ sufficiently close to $T_*$, we can assume $t_0$ itself is sufficiently close to $T_*$ and set $E_0(t)$ in \eqref{fE} to be
\beqs 
E_0(t)=M C_1^{-\delta}(T_*-t)^{\delta/\alpha}\text{ which has the form of \eqref{epower}. }
\eeqs 

\medskip
In Case 2, by \eqref{yplus} and \eqref{G2}, we have, as  $t<T_*$  close to $T_*$,
\beqs
|f(t)|\le \frac{M|y(t)|^{1+\alpha}}{\left(\ln(C_1(T_*-t)^{-1/\alpha}\right)^p}
\le \frac{M|y(t)|^{1+\alpha}}{\left(\frac12\ln((T_*-t)^{-1/\alpha}\right)^p}
= \frac{(2\alpha)^pM |y(t)|^{1+\alpha}}{\left |\ln(T_*-t)\right|^p}. 
\eeqs
We, again, can  assume $t_0$ is sufficiently close to $T_*$ and set $E_0(t)$ in \eqref{fE} to be
\beqs
E_0(t)=\frac{(2\alpha)^pM}{\left | \ln(T_*-t) \right |^p} \text{ which has the form of \eqref{elog}. } 
\eeqs

We apply  Theorem \ref{eg1} for Case 1, and Theorem \ref{eg2} for  Case 2. 
Then there exists an eigenvector $\xi_*$  of $A$ associated with an eigenvalue $\Lambda$ such that \eqref{xiHA} holds together with estimates  \eqref{sc1}  and \eqref{sc2} in the respective cases.
It is obvious that \eqref{mainest1} and \eqref{mainest2} follow \eqref{sc1}   and \eqref{sc2}, respectively.
\end{proof}

\begin{example}\label{bioapp}
 Consider the following simplistic model for populations of two species which can cooperate and compete with each others at the same time:
 \beq\label{biosys}
\left\{ 
\begin{aligned}y_1'&=k_1(y_1+y_2)^\alpha y_1-(a_1 + b_1 y_2^{\beta_1})y_1,\\
  y_2'&=k_2(y_1+y_2)^\alpha y_2-(a_2 + b_2 y_1^{\beta_2})y_2,
  \end{aligned}
  \right.
 \eeq 
 where $y_i(t)\ge 0$, for $i=1,2$, is the population of the $i$-th species. Above, $\alpha>0$, $k_i>0$, $a_i>0$, $b_i\ge 0$, $\beta_i\ge 0$ are constants. (It is  clear that one can ignore $\beta_i$ whenever $b_i=0$.)
For $i=1,2$, the term $k_i(y_1+y_2)^\alpha$ describes a rapid birth rate as in the Malthus model but can be curbed down when $\alpha<1$. 
 The dependence on the total population $N=y_1+y_2$ means that the two species help each others thrive. (See also \cite{KaKa2016,Karev2005} for models of inhomogeneous populations with a similar dependence.) The coefficient $a_i$ is the natural death rate while $b_i y_{3-i}^{\beta_i}$ is the death rate caused by the competition with the other species. 

 Assume $\alpha>\beta_1$ and $\alpha>\beta_2$, which can roughly be interpreted as the cooperation is weightier than the competition. 
Then we obtain the equation \eqref{yGeq} with $n=2$, $y=(y_1,y_2)$ and
 \beq\label{biofun}
 H(y)=(|y_1|+|y_2|)^\alpha,\ 
 A={\rm diag}[k_1,k_2],\ 
 G(t,y)=\begin{pmatrix}-(a_1 + b_1 |y_2|^{\beta_1})y_1\\ -(a_2 + b_2 |y_1|^{\beta_2})y_2
            \end{pmatrix}.
 \eeq
Note that the function $G$ satisfies condition \eqref{G1}. 

\medskip\noindent\textit{Existence of blow-up solutions.} Suppose $|y(0)|$ is sufficiently large and $y_1(0)\ge 0$ and $y_2(0)\ge 0$.
By  Theorem \ref{blowthm}, there exists a solution $y(t)\in \R^2$ of \eqref{yGeq} and \eqref{biofun}  for $t\in[0,T_*)$ which blows up at the finite time $T_*$.
For $i=1,2$, let 
$$g_i(t)=k_i(|y_1(t)|+|y_2(t)|)^\alpha -(a_i + b_i |y_{3-i}(t)|^{\beta_i}).$$
 We rewrite the system \eqref{yGeq} and \eqref{biofun} as
\beq
y_i'=g_i(t)y_i,\quad i=1,2.
\eeq
Hence, $y_i(t)=y_i(0)e^{\int_0^t g_i(\tau)\d \tau}$ for $t\in[0,T_*)$. This implies $y_1(t)\ge 0$ and $y_2(t)\ge 0$ for all $t\in[0,T_*)$, hence, $y(t)$ is a blow-up solution of \eqref{biosys} now.

\medskip\noindent\textit{Asymptotic behavior near the blow-up time.}
We can apply Theorem \ref{typthm}, Case 1. Therefore, any solution of \eqref{biosys} that blows up at a finite time $T_*$ behaves, when $t\to T_*^-$,  like
either 
\beq\label{beh1}
c (T_*-t)^{-1/\alpha} e_j \text{ for some $j\in\{1,2\}$ and number $c > 0$,  when $k_1\ne k_2$,}
\eeq
or 
\beq\label{beh2}
(T_*-t)^{-1/\alpha}\begin{pmatrix}a\\ b \end{pmatrix} \text{ for some numbers $a,b\ge 0$ with $a^2+b^2>0$,  when $k_1=k_2=k$.}
\eeq
Above $e_1=(1,0)$ and $e_2=(0,1)$. In fact, thanks to \eqref{xiHA}, we have $\alpha k_j c^\alpha=1$ in the case \eqref{beh1}, or $\alpha k (a+b)^\alpha=1$ in the case \eqref{beh2}. 
Consequently, when $k_1\ne k_2$, it follows from \eqref{beh1} that one species' population will go to infinity much faster than the other's. Whether the remaining species goes extinct requires the second approximation or even an asymptotic expansion of the solution.
\end{example}

\begin{remark}\label{comparelit}
 We compare our methods and results with the well-known ones used for parabolic equations. We refer to those in the established text \cite{GalakBook1995} on the subject.
 (See also a more recent book \cite{BeiHu2011} with the same approach and similar results.)
 
 \begin{enumerate}[label=\rnum]
  \item  Our techniques are developed genuinely for \textit{systems}. Those in \cite{GalakBook1995} deal with scalar equations albeit PDEs and can use of the Maximum and Comparison Principles that we do not have.
  \item Our approach yields the precise asymptotic behavior near the blow-up time for \textit{any} blow-up solution.
The approach in \cite{GalakBook1995} finds the blow-up profiles for the solutions in a certain class, called self-similar solutions,  and establishes their stability, i.e. studies their perturbed solutions. Thus, it does not deal with all blow-up solutions. 

\item Since our lower order term is very general, the self-similar solutions, in general, do not exist. Moreover, unlike \cite{GalakBook1995}, we do not use the scaling arguments.

\item Our current results do not cover those in \cite{GalakBook1995,BeiHu2011} and vice versa. It seems that the two approaches complement each others.
 \end{enumerate}
\end{remark}

\begin{remark}\label{findappn}
 The following final remarks are in order.
 
\begin{enumerate}[label=\rnum]
  \item The new classes of functions $G(t,x)$ in \eqref{yGeq} and $f(t)$ in \eqref{mainode},  and the proofs in this paper can be adapted to improve the results in previous work \cite{H7,H8}. 
  
  \item Our results provide details to the blow-up phenomena for some dynamical systems forward in time. 
    For dissipative systems even with nonlinear dissipative terms, although the blow-ups do not occur forward in time, the results in this paper can still be used to analyze the solutions with blow-ups backward in time.  These solutions \text{do not lie in} the global attractor. Therefore, this study will complement the approach that focuses on the solution \textit{on} the global attractor. Because of that, it will provide some clearer understanding of the complex dynamics.
    
  \item Although our motivation is a theory for PDEs, we are not able to achieve it yet and are settled for an ODE theory in this paper. At least, it verifies that the direction is promising. It may have a number of applications to ODE models of some population dynamics such as the one in Example \ref{bioapp}.
  Connecting to the PDEs from the ODE point of view, we hope to apply this theory to some finite-dimensional approximations of the PDEs.
  
  \item Improvements are certainly needed to strengthen the results' applicability, especially to PDEs.
\end{enumerate}
\end{remark}

\medskip

\textit{Acknowledgment.} The author would like to thank the anonymous referee for valuable suggestions that help improve the paper's presentation.

\medskip
\noindent\textbf{Data availability.} 
No new data were created or analyzed in this study.

\medskip
\noindent\textbf{Funding.} No funds were received for conducting this study. 

\medskip
\noindent\textbf{Conflict of interest.}
There are no conflicts of interests.

\bibliography{paperbaseall}{}

@preamble{"\def\cprime{$'$}"}

@book {BeiHu2011,
    AUTHOR = {Hu, Bei},
     TITLE = {Blow-up theories for semilinear parabolic equations},
    SERIES = {Lecture Notes in Mathematics},
    VOLUME = {2018},
 PUBLISHER = {Springer, Heidelberg},
      YEAR = {2011},
     PAGES = {x+125},
      ISBN = {978-3-642-18459-8},
   MRCLASS = {35B44 (35-02 35B50 35C06 35K58 35K91 35Q55)},
  MRNUMBER = {2796831},
MRREVIEWER = {Hongwei\ Chen},
       DOI = {10.1007/978-3-642-18460-4},
       URL = {https://doi.org/10.1007/978-3-642-18460-4},
}

@book{GalakBook1995,
url = {https://doi.org/10.1515/9783110889864},
title = {Blow-Up in Quasilinear Parabolic Equations},
author = {A. A. Samarskii and Victor A. Galaktionov and Sergey P. Kurdyumov and A. P. Mikhailov},
publisher = {De Gruyter},
address = {Berlin, New York},
doi = {doi:10.1515/9783110889864},
isbn = {9783110889864},
year = {1995},
lastchecked = {2025-09-01}
}

@article{Karev2005,
author = {Karev, George P.},
title = {Dynamics of inhomogeneous populations and global demography models},
journal = {Journal of Biological Systems},
volume = {13},
number = {01},
pages = {83-104},
year = {2005},
doi = {10.1142/S0218339005001410},
URL = {https://doi.org/10.1142/S0218339005001410},
eprint = {https://doi.org/10.1142/S0218339005001410}
}

@article {KaKa2016,
    AUTHOR = {G P Karev and I Kareva},
     TITLE = {Mathematical Modeling of Extinction of Inhomogeneous Populations},
   JOURNAL = {Bull Math Biol.},
  FJOURNAL = {Bulletin of Mathematical Biology},
    VOLUME = {78},
      YEAR = {2016},
    NUMBER = {4},
     PAGES = {834--858},
      ISSN = {},
     CODEN = {},
   MRCLASS = {},
  MRNUMBER = {},
MRREVIEWER = {},
       DOI = {10.1007/s11538-016-0166-0},
       URL = {},
}

@book {Hartman1964,
    AUTHOR = {Hartman, Philip},
     TITLE = {Ordinary differential equations},
 PUBLISHER = {John Wiley \& Sons, Inc., New York-London-Sydney},
      YEAR = {1964},
     PAGES = {xiv+612},
   MRCLASS = {34.00},
  MRNUMBER = {0171038},
MRREVIEWER = {H. A. Antosiewicz},
}

@article{Ghidaglia1986b,
	author = "Jean-Michel Ghidaglia",
	doi = "10.1016/0362-546X(86)90037-4",
	fjournal = "Nonlinear Analysis. Theory, Methods \& Applications. An International Multidisciplinary Journal",
	issn = "0362-546X",
	journal = "Nonlinear Anal.",
	mrclass = "34G20 (35K22 47H15)",
	mrnumber = "851146",
	mrreviewer = "Andrzej Hajnosz",
	number = "8",
	pages = "777--790",
	title = "Some backward uniqueness results",
	url = "https://doi.org/10.1016/0362-546X(86)90037-4",
	volume = "10",
	year = "1986"
}

@article{Ghidaglia1986a,
	author = "Jean-Michel Ghidaglia",
	doi = "10.1016/0022-0396(86)90121-X",
	fjournal = "Journal of Differential Equations",
	issn = "0022-0396",
	journal = "J. Differential Equations",
	mrclass = "35Q10 (35B40)",
	mrnumber = "823404",
	number = "2",
	pages = "268--294",
	title = "Long time behaviour of solutions of abstract inequalities: applications to thermohydraulic and magnetohydrodynamic equations",
	url = "https://doi.org/10.1016/0022-0396(86)90121-X",
	volume = "61",
	year = "1986"
}

@article{Shi2000,
	author = "Y. Shi",
	fjournal = "Communications in Partial Differential Equations",
	issn = "0360-5302",
	journal = "Comm. Partial Differential Equations",
	mrclass = "35L70 (34G20 35C20 47D06)",
	mrnumber = "1789928",
	mrreviewer = "Alfredo Marzocchi",
	number = "11-12",
	pages = "2287--2331",
	title = "A {F}oias-{S}aut type of expansion for dissipative wave equations",
	url = "https://doi.org/10.1080/03605300008821585",
	volume = "25",
	year = "2000"
}

@book{JHale78,
	address = "Huntington, N.Y.",
	author = "Jack K. Hale",
	edition = "Second",
	isbn = "0-89874-011-8",
	mrclass = "34-02 (58Fxx)",
	mrnumber = "587488 (82e:34001)",
	mrreviewer = "Charles C. Conley",
	pages = "xvi+361",
	publisher = "Robert E. Krieger Publishing Co. Inc.",
	title = "{Ordinary differential equations}",
	year = "1980"
}

@article{Minea,
	author = "Gheorghe Minea",
	coden = "JDDEEH",
	doi = "10.1023/A:1022696614020",
	fjournal = "Journal of Dynamics and Differential Equations",
	issn = "1040-7294",
	journal = "J. Dynam. Differential Equations",
	mrclass = "34C20 (37G05)",
	mrnumber = "MR1607538 (2000b:34058)",
	number = "1",
	pages = "189--207",
	title = "{Investigation of the {F}oias-{S}aut normalization in the finite-dimensional case}",
	url = "http://dx.doi.org/10.1023/A:1022696614020",
	volume = "10",
	year = "1998"
}

@article{FS84a,
	author = "C. Foias and J.-C. Saut",
	coden = "IUMJAB",
	fjournal = "Indiana University Mathematics Journal",
	issn = "0022-2518",
	journal = "Indiana Univ. Math. J.",
	mrclass = "35Q10 (35B40 76D05)",
	mrnumber = "MR740960 (85f:35164)",
	mrreviewer = "Charles J. Amick",
	number = "3",
	pages = "459--477",
	title = "{Asymptotic behavior, as {$t\rightarrow +\infty $}, of solutions of {N}avier-{S}tokes equations and nonlinear spectral manifolds}",
	volume = "33",
	year = "1984"
}

@article{FS84b,
	author = "C. Foias and J.-C. Saut",
	coden = "IUMJAB",
	fjournal = "Indiana University Mathematics Journal",
	issn = "0022-2518",
	journal = "Indiana Univ. Math. J.",
	mrclass = "35Q10 (76D05)",
	mrnumber = "MR763949 (86b:35168)",
	mrreviewer = "Charles J. Amick",
	number = "6",
	pages = "911--926",
	title = "{On the smoothness of the nonlinear spectral manifolds associated to the {N}avier-{S}tokes equations}",
	volume = "33",
	year = "1984"
}

@article{FGuillope86,
	author = "C. Foias and C. Guillop{\'e}",
	coden = "JDEQAK",
	fjournal = "Journal of Differential Equations",
	issn = "0022-0396",
	journal = "J. Differential Equations",
	mrclass = "35Q10 (35B35 35B40 76D05)",
	mrnumber = "MR818863 (87g:35187)",
	mrreviewer = "Howard Swann",
	number = "1",
	pages = "128--148",
	title = "{On the behavior of the solutions of the {N}avier-{S}tokes equations lying on invariant manifolds}",
	volume = "61",
	year = "1986"
}

@incollection{FS86,
	address = "Providence, RI",
	author = "C. Foias and J.-C. Saut",
	booktitle = "{Nonlinear functional analysis and its applications, Part 1 (Berkeley, Calif., 1983)}",
	mrclass = "35Q10 (76D05)",
	mrnumber = "MR843577 (88f:35122)",
	mrreviewer = "C. Bardos",
	pages = "439--448",
	publisher = "Amer. Math. Soc.",
	series = "{Proc. Sympos. Pure Math.}",
	title = "{Nonlinear spectral manifolds for the {N}avier-{S}tokes equations}",
	volume = "45",
	year = "1986"
}

@article{FS87,
	author = "C. Foias and J.-C. Saut",
	fjournal = "Annales de l'Institut Henri Poincar{\'e}. Analyse Non Lin{\'e}aire",
	issn = "0294-1449",
	journal = "Ann. Inst. H. Poincar{\'e} Anal. Non Lin{\'e}aire",
	mrclass = "35Q10 (35B40 35Q20 76D05)",
	mrnumber = "MR877990 (88d:35158)",
	mrreviewer = "Michael Wiegner",
	number = "1",
	pages = "1--47",
	title = "{Linearization and normal form of the {N}avier-{S}tokes equations with potential forces}",
	volume = "4",
	year = "1987"
}

@article{FS91,
	author = "C. Foias and J.-C. Saut",
	coden = "IUMJAB",
	fjournal = "Indiana University Mathematics Journal",
	issn = "0022-2518",
	journal = "Indiana Univ. Math. J.",
	mrclass = "35Q30 (35C20 76D05)",
	mrnumber = "MR1101233 (92g:35175)",
	mrreviewer = "Song Mu Zheng",
	number = "1",
	pages = "305--320",
	title = "{Asymptotic integration of {N}avier-{S}tokes equations with potential forces. {I}}",
	volume = "40",
	year = "1991"
}

@article {CaH3,
    AUTHOR = {Cao, Dat and Hoang, Luan},
     TITLE = {Asymptotic expansions with exponential, power, and logarithmic
              functions for non-autonomous nonlinear differential equations},
   JOURNAL = {J. Evol. Equ.},
  FJOURNAL = {Journal of Evolution Equations},
    VOLUME = {21},
      YEAR = {2021},
    NUMBER = {2},
     PAGES = {1179--1225},
      ISSN = {1424-3199},
   MRCLASS = {34E05 (34E10 41A60)},
  MRNUMBER = {4278393},
       DOI = {10.1007/s00028-020-00622-w},
       URL = {https://doi.org/10.1007/s00028-020-00622-w},
}

@article {CaHK1,
    AUTHOR = {Cao, Dat and Hoang, Luan and Kieu, Thinh},
     TITLE = {Infinite series asymptotic expansions for decaying solutions
              of dissipative differential equations with non-smooth
              nonlinearity},
   JOURNAL = {Qual. Theory Dyn. Syst.},
  FJOURNAL = {Qualitative Theory of Dynamical Systems},
    VOLUME = {20},
      YEAR = {2021},
    NUMBER = {3},
     PAGES = {Paper No. 62, 38 pp},
      ISSN = {1575-5460},
   MRCLASS = {35Q30 (34E05 41A60 76D05)},
  MRNUMBER = {4284023},
       DOI = {10.1007/s12346-021-00502-9},
       URL = {https://doi.org/10.1007/s12346-021-00502-9},
}

@article {CaH2,
    AUTHOR = {Cao, Dat and Hoang, Luan},
     TITLE = {Asymptotic expansions in a general system of decaying
              functions for solutions of the {N}avier-{S}tokes equations},
   JOURNAL = {Ann. Mat. Pura Appl. (4)},
  FJOURNAL = {Annali di Matematica Pura ed Applicata. Series IV},
    VOLUME = {199},
      YEAR = {2020},
    NUMBER = {3},
     PAGES = {1023--1072},
      ISSN = {0373-3114},
   MRCLASS = {76D05 (35B40 35C20 35Q30 37L05 76M45)},
  MRNUMBER = {4102801},
       DOI = {10.1007/s10231-019-00911-3},
       URL = {https://doi.org/10.1007/s10231-019-00911-3},
}

@article {CaH1,
    AUTHOR = {Cao, Dat and Hoang, Luan},
     TITLE = {Long-time asymptotic expansions for {N}avier-{S}tokes
              equations with power-decaying forces},
   JOURNAL = {Proc. Roy. Soc. Edinburgh Sect. A},
  FJOURNAL = {Proceedings of the Royal Society of Edinburgh. Section A.
              Mathematics},
    VOLUME = {150},
      YEAR = {2020},
    NUMBER = {2},
     PAGES = {569--606},
      ISSN = {0308-2105},
   MRCLASS = {35Q30 (35B40 35C20 76D05)},
  MRNUMBER = {4080452},
       DOI = {10.1017/prm.2018.154},
       URL = {https://doi.org/10.1017/prm.2018.154},
}

@article{HM2,
	author = "Luan T. Hoang and Vincent R. Martinez",
	doi = "10.1016/j.jmaa.2018.01.065",
	fjournal = "Journal of Mathematical Analysis and Applications",
	journal = "J. Math. Anal. Appl.",
	number = "1",
	pages = "84--113",
	title = "Asymptotic expansion for solutions of the {N}avier-{S}tokes equations with non-potential body forces",
	url = "https://doi.org/10.1016/j.jmaa.2018.01.065",
	volume = "462",
	year = "2018",
}

@article {HTi1,
    AUTHOR = {Hoang, Luan T. and Titi, Edriss S.},
     TITLE = {Asymptotic expansions in time for rotating incompressible
              viscous fluids},
   JOURNAL = {Ann. Inst. H. Poincar\'{e} Anal. Non Lin\'{e}aire},
  FJOURNAL = {Annales de l'Institut Henri Poincar\'{e}. Analyse Non Lin\'{e}aire},
    VOLUME = {38},
      YEAR = {2021},
    NUMBER = {1},
     PAGES = {109--137},
      ISSN = {0294-1449},
   MRCLASS = {35Q30 (35C20 76D05 76E07 76U05)},
  MRNUMBER = {4200479},
       DOI = {10.1016/j.anihpc.2020.06.005},
       URL = {https://doi.org/10.1016/j.anihpc.2020.06.005},
}

@article{HM1,
	author = "Luan T. Hoang and Vincent R. Martinez",
	doi = "10.3233/ASY-171429",
	fjournal = "Asymptotic Analysis",
	issn = "",
	journal = "Asymptot. Anal.",
	mrclass = "",
	mrnumber = "",
	number = "3--4",
	pages = "167--190",
	title = "Asymptotic expansion in {G}evrey spaces for solutions of {N}avier-{S}tokes equations",
	url = "http://dx.doi.org/10.3233/ASY-171429",
	volume = "104",
	year = "2017",
}

@article{FHS1,
	author = "Ciprian Foias and Luan Hoang and Jean-Claude Saut",
	coden = "JFUAAW",
	doi = "10.1016/j.jfa.2011.02.005",
	fjournal = "Journal of Functional Analysis",
	issn = "0022-1236",
	journal = "J. Funct. Anal.",
	mrclass = "37Nxx (35Q30 37Lxx 76Dxx)",
	mrnumber = "2774064",
	number = "10",
	pages = "3007--3035",
	title = "{Asymptotic integration of {N}avier-{S}tokes equations with potential forces. {II}. {A}n explicit {P}oincar{\'e}-{D}ulac normal form}",
	url = "http://dx.doi.org/10.1016/j.jfa.2011.02.005",
	volume = "260",
	year = "2011",
}

@article{FHOZ2,
	author = "Ciprian Foias and Luan Hoang and Eric Olson and Mohammed Ziane",
	doi = "10.1016/j.anihpc.2008.09.003",
	fjournal = "Annales de l'Institut Henri Poincar{\'e}. Analyse Non Lin{\'e}aire",
	issn = "0294-1449",
	journal = "Ann. Inst. H. Poincar{\'e} Anal. Non Lin{\'e}aire",
	mrclass = "35Qxx (76D05)",
	mrnumber = "MR2566704",
	number = "5",
	pages = "1635--1673",
	title = "{The normal form of the {N}avier-{S}tokes equations in suitable normed spaces}",
	url = "http://dx.doi.org/10.1016/j.anihpc.2008.09.003",
	volume = "26",
	year = "2009",
}

@article{FHN2,
	author = "Ciprian Foias and Luan Hoang and Basil Nicolaenko",
	coden = "CMPHAY",
	doi = "10.1007/s00220-009-0827-z",
	fjournal = "Communications in Mathematical Physics",
	issn = "0010-3616",
	journal = "Comm. Math. Phys.",
	mrclass = "35Q30 (35B10 35B40 35C20 76D03 76D05 76F02 76M35)",
	mrnumber = "MR2525635 (2010f:35281)",
	mrreviewer = "Peter E. Kloeden",
	number = "2",
	pages = "679--717",
	title = "{On the helicity in 3{D}-periodic {N}avier-{S}tokes equations. {II}. {T}he statistical case}",
	url = "http://dx.doi.org/10.1007/s00220-009-0827-z",
	volume = "290",
	year = "2009",
}

@article{FHN1,
	author = "Ciprian Foias and Luan Hoang and Basil Nicolaenko",
	doi = "10.1112/plms/pdl003",
	fjournal = "Proceedings of the London Mathematical Society. Third Series",
	issn = "0024-6115",
	journal = "Proc. Lond. Math. Soc. (3)",
	mrclass = "35Q30 (35B10 35B40 35C20 37L99 76D05)",
	mrnumber = "MR2293465 (2007k:35357)",
	mrreviewer = "Peter E. Kloeden",
	number = "1",
	pages = "53--90",
	title = "{On the helicity in 3{D}-periodic {N}avier-{S}tokes equations. {I}. {T}he non-statistical case}",
	url = "http://dx.doi.org/10.1112/plms/pdl003",
	volume = "94",
	year = "2007",
}

@article{FHOZ1,
	author = "Ciprian Foias and Luan Hoang and Eric Olson and Mohammed Ziane",
	coden = "IUMJAB",
	doi = "10.1512/iumj.2006.55.2830",
	fjournal = "Indiana University Mathematics Journal",
	issn = "0022-2518",
	journal = "Indiana Univ. Math. J.",
	mrclass = "35Q30 (35C20 76D03 76D05)",
	mrnumber = "MR2225448 (2008d:35151)",
	mrreviewer = "B. Szafirski",
	number = "2",
	pages = "631--686",
	title = "{On the solutions to the normal form of the {N}avier-{S}tokes equations}",
	url = "http://dx.doi.org/10.1512/iumj.2006.55.2830",
	volume = "55",
	year = "2006",
}

@article{H5,
	author = "Luan Hoang",
	journal = "Ann. Sc. Norm. Super. Pisa Cl. Sci.",
	fjournal = "The Annali della Scuola Normale di Pisa - Classe di Scienze",
	mrclass = "",
	mrnumber = "",
	doi = "10.2422/2036-2145.202109\_004",
	number = "1",
	pages = "311--370",
	title = "Asymptotic expansions about infinity for solutions of nonlinear differential equations with coherently decaying forcing functions",
	volume = "XXV",
	year = "2024",
}

@article{H6,
title = {The {N}avier--{S}tokes equations with body forces decaying coherently in time},
JOURNAL = {J. Math. Anal. Appl.},
FJOURNAL = {Journal of Mathematical Analysis and Applications},
volume = {531},
number = {2, Part 1},
pages = {Paper No. 127863, 39 pp},
year = {2024},
issn = {0022-247X},
doi = {https://doi.org/10.1016/j.jmaa.2023.127863},
url = {https://www.sciencedirect.com/science/article/pii/S0022247X23008661},
author = {Luan Hoang},
}

@article{H7,
author = {Luan Hoang},
title = {Long-time behaviour of solutions of superlinear systems of differential equations},
journal = {Dynamical Systems},
volume = {39},
number = {1},
pages = {79--107},
year = {2024},
publisher = {Taylor & Francis},
doi = {10.1080/14689367.2023.2234845},
URL = {https://doi.org/10.1080/14689367.2023.2234845},
eprint = {https://doi.org/10.1080/14689367.2023.2234845}
}

@article{H8,
	author = "Luan Hoang",
	fjournal = "Electronic Journal of Differential Equations",
	journal = "Electron.  J. Differential Equations",
	mrclass = "",
	mrnumber = "",
	note = "",
	number = "08",
	pages = "1--26",
	title = "Behavior near the extinction time for systems of differential equations with sublinear dissipation terms",
	volume = "2025",
	year = "2025",
	doi="10.58997/ejde.2025.08",
}
\bibliographystyle{plain}

 \end{document}